\documentclass{eptcs}
\providecommand{\event}{ACT2024}

\usepackage{doi}
\usepackage[backend=biber, backref=true, maxbibnames = 10, style = alphabetic]{biblatex}

\AtBeginBibliography{\small}

\usepackage{iftex}
\usepackage{amssymb}
\usepackage{amsmath}
\usepackage{mathtools}
\usepackage{amsthm}
\usepackage{amsfonts}
\usepackage{tikz-cd} 
\usepackage{graphicx} 
\usepackage{stmaryrd} 
\usepackage[bibliography=common]{apxproof} 
\usepackage{multicol}

\DeclareMathAlphabet{\dutchcal}{U}{dutchcal}{m}{n}
\SetMathAlphabet{\dutchcal}{bold}{U}{dutchcal}{b}{n}
\DeclareMathAlphabet{\dutchbcal} {U}{dutchcal}{b}{n}

\newtheoremstyle{theoremdd}
{\topsep}{\topsep}{\upshape}{0pt}{\bfseries}{.}{ }{\thmname{#1}\thmnumber{ #2}\thmnote{ (#3)}}

\theoremstyle{definition}
\newtheoremrep{Th}{Theorem}[section]
\newtheoremrep{Lemma}[Th]{Lemma}
\newtheoremrep{Cor}[Th]{Corollary}
\newtheoremrep{Prop}[Th]{Proposition}
\newtheoremrep{Def}[Th]{Definition}

\newtheorem{Rem}[Th]{Remark}
\newtheorem{Ex}[Th]{Example}
\newtheorem{Not}[Th]{Notation}

\newcommand{\cat}{\dutchcal} 
\newcommand{\ncat}{\mathbf} 
\DeclareMathSymbol{\hyp}{\mathord}{AMSa}{"39}
\newcommand{\Hom}{\text{Hom}}

\newcommand{\op}{\text{op}}

\newcommand{\Curr}{\text{Curr}}

\newcommand{\llp}{\llparenthesis}
\newcommand{\rrp}{\rrparenthesis}
\newcommand{\semp}[1]{\llp #1 \rrp}

\newcommand{\can}{\text{Can}}
\newcommand{\Fun}{\text{Fun}}
\newcommand{\Sort}{\text{Sort}}
\newcommand{\Term}{\text{Term}}

\newcommand{\Path}{\text{Path}}
\newcommand{\Free}{\text{Free}}

\newcommand{\Gen}{\text{Gen}}

\newcommand{\inst}[1]{\cat{#1}}
\newcommand{\CPath}{\text{CPath}}

\newcommand{\ccomp}{\circledast}
\newcommand{\bbprof}{\mathbb{P}\ncat{rof}}
\newcommand{\bbcurr}{\mathbb{C}\ncat{urr}}
\newcommand{\pro}{\slashedrightarrow}

\makeatletter
\def\arrowfillmod@#1#2#3#4#5{%
  $\m@th\thickmuskip0mu\medmuskip\thickmuskip\thinmuskip\thickmuskip
   \relax#5#1\mkern-7mu%
   \cleaders\hbox{$#5\mkern-2mu#2\mkern-2mu$}\hfill
   \mathclap{#3}\mathclap{#2}%
   \cleaders\hbox{$#5\mkern-2mu#2\mkern-2mu$}\hfill
   \mkern-7mu#4$%
}

\def\rightslashedarrowfill@{%
  \arrowfillmod@\relbar\relbar\mapstochar\rightarrow}
\newcommand\xslashedrightarrow[2][]{%
  \ext@arrow 0055{\rightslashedarrowfill@}{#1}{#2}}
\def\slashedrightarrow{\xslashedrightarrow{}}

\def\pgf@stroke@inner@line{%
  \let\pgf@temp@save=\pgf@strokecolor@global
  \pgfsys@beginscope%
  {%
    \pgfsys@roundcap
    \pgfsys@setlinewidth{\pgfinnerlinewidth}%
    \pgfsetstrokecolor{\pgfinnerstrokecolor}%
    \pgfsyssoftpath@invokecurrentpath%
    \pgfsys@stroke%
  }%
  \pgfsys@endscope%
  \global\let\pgf@strokecolor@global=\pgf@temp@save
}
\makeatother

\usetikzlibrary{calc}
\usetikzlibrary{decorations.pathmorphing}

\tikzset{curve/.style={settings={#1},to path={(\tikztostart)
    .. controls ($(\tikztostart)!\pv{pos}!(\tikztotarget)!\pv{height}!270:(\tikztotarget)$)
    and ($(\tikztostart)!1-\pv{pos}!(\tikztotarget)!\pv{height}!270:(\tikztotarget)$)
    .. (\tikztotarget)\tikztonodes}},
    settings/.code={\tikzset{quiver/.cd,#1}
        \def\pv##1{\pgfkeysvalueof{/tikz/quiver/##1}}},
    quiver/.cd,pos/.initial=0.35,height/.initial=0}

\tikzset{tail reversed/.code={\pgfsetarrowsstart{tikzcd to}}}
\tikzset{2tail/.code={\pgfsetarrowsstart{Implies[reversed]}}}
\tikzset{2tail reversed/.code={\pgfsetarrowsstart{Implies}}}

\ifpdf
  \usepackage{underscore}         
  \usepackage[T1]{fontenc}        
\else
  \usepackage{breakurl}           
\fi

\title{Presenting Profunctors}
\author{Gabriel Goren-Roig
\institute{Universidad de Buenos Aires \\ CONICET}
\email{ggoren@dc.uba.ar}
\and
Joshua Meyers
\institute{Conexus AI}
\email{meygerjos@protonmail.com}
\and
Emilio Minichiello
\institute{CUNY CityTech}
\email{eminichiello67@gmail.com}
}
\def\titlerunning{Presenting Profunctors}
\def\authorrunning{G. Goren-Roig, J. Meyers, \& E. Minichiello}
\begin{document}
\maketitle

\begin{abstract}
Motivated by problems in categorical database theory, we introduce and compare two notions of presentation for profunctors, uncurried and curried, which arise intuitively from thinking of profunctors either as functors $\cat{C}^\op \times \cat{D} \to \ncat{Set}$ or $\cat{C}^\op \to \ncat{Set}^{\cat{D}}$.
Although the Cartesian closure of $\ncat{Cat}$ means these two perspectives can be used interchangeably at the semantic level, a surprising amount of subtlety is revealed when looking through the lens of syntax. Indeed, we prove that finite uncurried presentations are strictly more expressive than finite curried presentations, hence the two notions do not induce the same class of finitely presentable profunctors. Moreover, an explicit construction for the composite of two curried presentations shows that the class of finitely curried presentable profunctors is closed under composition, in contrast with the larger class of finitely uncurried presentable profunctors, which is not. This shows that curried profunctor presentations are more appropriate for computational tasks that use profunctor composition. We package our results on curried profunctor presentations into a double equivalence from a syntactic double category into the double category of profunctors. Finally, we study the relationship between curried and uncurried presentations, leading to the introduction of curryable presentations. These constitute a subcategory of uncurried presentations which is equivalent to the category of curried presentations, therefore acting as a bridge between the two syntactic choices.
\end{abstract}

\paragraph{Acknowledgements:} We thank the reviewers for their detailed comments and suggestions. We also thank David Spivak and Ryan Wisnesky for helpful conversations.

\section{Introduction}

In the most recent iteration of categorical database theory \cite{spivak2012functorial, spivak2015relational, patterson2022categorical, schultz2017algebraic, schultz2017integration}, profunctors are pervasive. Indeed, they are used to define database schemas with attributes as well as to define generalized queries on such schemas, which, fundamentally, can also implement data exchange operations. In this context, composition of profunctors plays a key role, since it allows us to define composite queries in a data-independent way, and it specializes to \emph{query coevaluation} \cite{schultz2017algebraic}, a natural categorical operation on database instances which may have uses in data integration.

In order to encode these entities as finite data and define algorithms on them that implement operations such as composition, we must move into the territory of finite presentations of these mathematical objects, and, as we will show, different perspectives on profunctors naturally suggest distinct, inequivalent notions of presentation.
In \cite{schultz2017integration}, the authors develop a syntactic theory of data integration based on the concept of \emph{uberflowers}, which can be thought of as a notion of presentation for a certain generalization of profunctors. Note that uberflowers are central to the database algorithms implemented in the open-source CQL tool available at \href{http://categoricaldata.net}{categoricaldata.net}.
When taking the particular case of profunctors, this syntactic notion reduces to what we call \emph{curried profunctor presentations}. These are natural when thinking of profunctors $\cat{P}: \cat{C} \pro \cat{D}$ as functors $\cat{P}: \cat{C}^\op \to \ncat{Set}^{\cat{D}}$.

However, existing work on uberflowers does not connect these syntactic entities with their semantics as (generalized) profunctors. On the contrary, although the semantic counterpart of the theory of uberflowers is developed in \cite{schultz2017algebraic}, there the authors take the alternative view of profunctors as functors of the form $\cat{P} : \cat{C}^\op \times \cat{D} \to \ncat{Set}$. From this perspective, they define a notion of presentation for their semantic version of queries, which they call \emph{bimodule presentation}, and which reduces in the case of profunctors to what we call \emph{uncurried profunctor presentations}.


Owing to the Cartesian closure of $\ncat{Cat}$, these two descriptions of a profunctor are interchangeable. However, we will show that the corresponding syntactic notions of curried and uncurried profunctor presentations are not, since they give rise to inequivalent notions of finite presentability, and hence \emph{a fortiori} uberflowers and bimodule presentations are not equivalent notions.

In Section \ref{section category presentations}, we discuss presentations of small categories.
In Section \ref{section profunctor presentations}, after having reviewed profunctors, we introduce our main objects of study: uncurried and curried profunctor presentations.
We begin our exposition with uncurried profunctor presentations (Definition \ref{def uncurried profunctor presentation})
and prove (Proposition \ref{prop uncurried prof don't admit finite composition}) that the class of finitely uncurried presentable profunctors is not closed under composition.
This already implies that bimodule presentations as defined in \cite{schultz2017algebraic} are not suitable for implementation. 

We then derive presentations for copresheaves on a category presentation $C$, which we call \emph{$C$-instance presentations}, as a particular case of uncurried profunctor presentations (Definition \ref{def instance presentation}). This is done in order to finally define curried presentations as indexed families of instance presentations and morphisms between them (Definition \ref{def curried profunctor presentation}). We define a composition operation for curried profunctor presentations and prove its correctness (Theorem \ref{th curried composition semantics}), which constitutes one of our main results. This immediately implies that the class of finitely curried presentable profunctors is closed under composition.
We conclude that curried presentations are indeed appropriate for computing profunctor composition.
Although the idea behind composition of curried profunctor presentations is present in \cite{schultz2017integration}, the precise mathematical formulation and the proof of correctness are novel to the best of our knowledge.

%
%
%



In Section \ref{section currying}, we turn towards the question of how exactly these two notions of profunctor presentation relate to each other. We begin by noting that there is an evident ``uncurrying'' functor from curried to uncurried presentations which preserves the semantics and sends finite presentations to finite presentations (Proposition \ref{prop curr to uncurr}). This, together with 
the fact that there are finitely uncurried presentable profunctors which are not finitely curried presentable (Lemma \ref{lem finite uncurried but not finite curried}),
lets us conclude that finite uncurried presentations are strictly more expressive than their curried counterparts. We then study the image of the uncurrying functor, and in doing so we arrive at the concept of \emph{curryable profunctor presentations}.
%
In our second main theorem (Theorem \ref{th curryable equiv to curried}), we prove that the resulting category of curryable profunctor presentations is equivalent to the category of curried presentations, and that this equivalence restricts to the respective full subcategories of finite presentations. This shows that we can indeed replace finite curried presentations by finite curryable presentations, which may be more convenient to work with.


Finally, in Section \ref{section double categories} we package our results on curried profunctor presentations into a double functor $\semp{-}: \mathbb{C}\ncat{urr} \to \mathbb{P}\ncat{rof}$ of double categories, and show that by appropriately quotienting the domain of this double functor, we obtain a double equivalence.

For clarity and space considerations, we defer most proofs of this paper to the Appendix.

%

\section{Category Presentations} \label{section category presentations}

In this section, we define category presentations together with their semantics via a functor $\semp{-}: \ncat{CatPr} \to \ncat{Cat}$. This functor factors through an appropriate quotient $\ncat{CatPr}_\approx$ of $\ncat{CatPr}$ and induces an equivalence $\ncat{CatPr}_\approx \simeq \ncat{Cat}$. We follow the terminology of algebraic theories (sorts, function symbols and signatures). Similar definitions can be found in \cite{spivak2014category}. 

\begin{Def}\label{def category presentation}
A \textbf{category signature} $\Sigma$ consists of a set $\Sort(\Sigma)$, whose elements we call sorts, and a set $\Fun(\Sigma)$, whose elements we call function symbols, along with two functions $s,t: \Fun(\Sigma) \to \Sort(\Sigma)$, called the source and target maps. We write $f: c_0 \to c_1$ to denote that $s(f) = c_0$ and $t(f) = c_1$. 

Given a category signature $\Sigma$ and $c, c' \in \Sort(\Sigma)$, a \textbf{path in $\Sigma$ from $c$ to $c'$}, denoted as $p: c \to c'$, consists of a (possibly empty) finite, ordered list of function symbols $p = (f_0, f_1, \dots, f_{n-1})$, such that if it is nonempty then $s(f_0) = c$, $t(f_{n-1}) = c'$ and $t(f_i) = s(f_{i+1})$ for all $0 \leq i < n-1$, and if it is empty we require $c = c'$. In the first case we write $p = f_0.f_1.\cdots.f_{n-1}$ while in the second case we write $p = 1_c$. Let $\Path(\Sigma)(c,c')$ denote the set of paths in $\Sigma$ from $c$ to $c'$ and $\Path(\Sigma)$ denote the set of all paths.

Given a category signature $\Sigma$, an \textbf{equation} over $\Sigma$ consists of a pair $(p_0, p_1)$ where $p_0,p_1: c \to c'$ are paths in $\Sigma$ with the same source and target. We denote such an equation by $p_0 = p_1: c \to c'$.

A \textbf{category presentation} is a pair $C = (C_\Sigma, C_E)$ of a category signature $C_\Sigma$ and a set $C_E$ of equations over $C_\Sigma$. We say that $C$ is a \textbf{finite category presentation} if $\Sort(C_\Sigma)$, $\Fun(C_\Sigma)$ and $C_E$ are all finite sets. 
\end{Def}

\begin{Not}
If $C=(C_\Sigma,C_E)$ is a category presentation, then we let $\Sort(C) \coloneqq \Sort(C_\Sigma)$, $\Fun(C) \coloneqq \Fun(C_\Sigma)$, and $\Path(C) \coloneqq \Path(C_\Sigma)$. We write $c \in C$ to mean that $c \in \Sort(C)$. Given $c,c' \in C$, we write $\Fun(C)(c,c')$ to denote the subset of function symbols $f$ of $C$ with $s(f) = c$ and $t(f) = c'$. We write $p_0 =_C p_1$ to mean that the equation $p_0 = p_1$ belongs to $C_E$. Since the equality symbol is taken up for equations in category presentations, given two paths $p$ and $q$ we use the notation $p \equiv q$ to indicate true (syntactic) equality. We still use $=$ for true equality of sorts and semantic objects.
\end{Not}

Equations in a category presentation constitute axioms from which equivalence of paths may be deduced by the rules of equational logic. This is captured by the following definition.

\begin{toappendix}
    We repeat here, for ease of reference, the definition of provable equality.
\end{toappendix}
\begin{Defrep}\label{def approx characterization}
Given a category presentation $C$, define a relation $\approx_C$ on paths of $C$ inductively\footnote{That is to say, there are no pairs $(p,q)$ in $\approx_C$ other than what can be derived in finitely many steps with the rules 1-6.} by the following inference rules:
\begin{enumerate}
    \item For any equation $p =_C q$, $p \approx_C q$.
    \item For any path $p$, $p\approx_C p$.
    \item If $p\approx_C q$, then $q\approx_C p$.
    \item If $p\approx_C q$ and $q\approx_C r$, then $p\approx_C r$.
    \item For any $C$-function symbol $f: c' \to c''$, if $p \approx_C q : c\to c'$, then $p.f \approx_C q.f$.
    \item For any $C$-function symbol $f: c \to c'$, if $p \approx_C q : c'\to c''$, then $f.p \approx_C f.q$.
\end{enumerate}

We call $\approx_C$ the \textbf{$C$-provable equality relation} on paths, and we say that $p_0$ and $p_1$ are $C$-provably equal if $p_0 \approx_C p_1$.
\end{Defrep}

\begin{toappendix}
\begin{Rem}\label{rem structural induction on inference rules}
Frequently we will use structural induction
on the rules 1-6 in Definition \ref{def approx characterization} to prove statements of the form ``for all paths $p$ and $q$, if $p\approx_C q$, then $P(p,q)$'' where $P(p,q)$ is some property of $p$ and $q$. To write such a proof, we simply need to show that for every inference rule, if the pairs of paths in the antecedent satisfy $P$, then the pair of paths in the conclusion satisfies $P$ as well. For example, for rule 6, we would show that for all parallel $C$-paths $p$ and $q$ and $C$-function symbols $f$ composable with $p$ and $q$, if $p \approx_C q$ and $P(p,q)$, then $P(f.p, f.q)$.

Notice that if the property $P(p,q)$ is itself an equivalence relation, then the verifications corresponding to rules 2-4 are satisfied automatically.
\end{Rem}
\end{toappendix}

We now define morphisms between category presentations. These map function symbols to paths in a way that sends equations (axioms) to provable equalities (equational theorems obtained from those axioms).

\begin{Def} \label{def morphism of category presentations}
Given category signatures $\Sigma$ and $\Sigma'$, a \textbf{morphism of category signatures} $F: \Sigma \to \Sigma'$ consists of functions $F_\Sort: \Sort(\Sigma) \to \Sort(\Sigma')$ and $F_\Fun: \Fun(\Sigma) \to \Path(\Sigma')$, such that if $f : c \to c'$ is a function symbol in $\Sigma$, then $F_\Fun(f): F_\Sort(c) \to F_\Sort(c')$. We refer to both $F_\Sort$ and $F_\Fun$ by $F$.
We extend $F$ to $\Path(\Sigma)$ by setting $F(f_0.f_1.\cdots.f_{n-1}) :\equiv F(f_0).F(f_1).\cdots.F(f_{n-1})$ and $F(1_c):\equiv 1_{F(c)}$.

Given category presentations $C$ and $C'$, a \textbf{morphism of category presentations} $F: C \to C'$ is a morphism $F: C_\Sigma \to C_\Sigma'$ of category signatures such that for each equation $p_0 =_C p_1$, we have $F(p_0) \approx_{C'} F(p_1)$. Let $\ncat{CatPr}$ denote the category of category presentations. Notice that the action of morphisms of category presentations on arbitrary paths is necessary for composition to be well defined.
\end{Def}

We obtain a category out of a category presentation $C$ by performing a quotient of the category freely generated by $C_\Sigma$ according to $C_E$. To that end we must first define congruences of categories (see \cite[Section 5.1]{borceux1994handbook} and \cite{bednarczyk1999generalized}).

\begin{Def}
Let $\cat{C}$ be a small category, and $\sim$ an equivalence relation on $\text{Mor}(\cat{C})$. We say that $\sim$ is a \textbf{congruence on} $\cat{C}$ if whenever $f \sim g$, $f$ and $g$ share the same domain and codomain; and furthermore, if $f_1 \sim f_2$, $g_1 \sim g_2$ and the composition $g_1 \circ f_1$ exists, then $g_1 \circ f_1 \sim g_2 \circ f_2$.  If $f$ is a morphism in $\cat{C}$, then we denote its $\sim$-equivalence class by $[f]$.

If $\sim$ is a congruence on $\cat{C}$, let $\cat{C}/{\sim}$ denote the small category whose objects are the same as $\cat{C}$, whose set of morphisms is $\text{Mor}(\cat{C})/{\sim}$, and where we define $1_c\coloneqq [1_c]$ and $[g]\circ [f]\coloneqq [g\circ f]$. We call it the \textbf{quotient category} of $\cat{C}$ by $\sim$ and denote it by $\cat{C}/{\sim}$.

If $R$ is a relation on $\text{Mor}(\cat{C})$ such that whenever $fRg$, $f$ and $g$ share the same domain and codomain, let $\sim_R$ be the smallest congruence that contains $R$.  Then we define $\cat{C}/R \coloneqq \cat{C}/{\sim_R}$, and we let $p_R : \cat{C} \to (\cat{C}/R)$ denote the canonical projection functor.
\end{Def}

\begin{Lemma} \label{lem quotient category universal property}
Let $\cat{C}$ be a small category and $R$ be a relation on $\text{Mor}(\cat{C})$ such that whenever $fRg$, $f$ and $g$ share the same domain and codomain.  Suppose that $F:\cat{C}\to\cat{D}$ is a functor with the property that if $f$ and $g$ are morphisms in $\cat{C}$ such that $fRg$, then $F(f) = F(g)$. Then there exists a unique functor $\widetilde{F}: (\cat{C}/R) \to \cat{D}$ such that $\widetilde{F} \circ p_R = F$.
\end{Lemma}

\begin{Lemmarep}\label{lem approx sim coincide}
Given a category presentation $C$, $\approx_C$ is the smallest congruence containing $C_E$. In other words, the relations $\sim_{C_E}$ and $\approx_C$ coincide.
\end{Lemmarep}

\begin{proof}
By rule 1, $\approx_C$ contains $C_E$, while by rules 2, 3, and 4, $\approx_C$ is an equivalence relation. To establish that $\approx_C$ is a congruence it remains to prove that $p\approx_C q$ implies $p$ and $q$ have the same domain and codomain (which follows straightforwardly by induction) and that $\approx_C$ is preserved by composition. To this end, we extend rules 5 and 6 to the case where $f$ is a path, this time by induction on $f$. Then whenever $f_1\approx_C f_2$ and $g_1\approx_C g_2$ and the composition $f_1.g_1$ exists, we have $f_1.g_1\approx_C f_1.g_2$ and $f_1.g_2\approx_C f_2.g_2$; then we conclude $f_1.g_1\approx f_2.g_2$ by rule 4. Now let $\approx$ be a congruence containing $C_E$. A straightforward application of induction proves that $p\approx_C q$ implies $p\approx q$. So $\approx_C$ is the smallest congruence containing $C_E$.
\end{proof}

\begin{Def} \label{def semantic category}
Given a category signature $\Sigma$, let $\Free(\Sigma)$ denote the small category whose objects are the sorts of $\Sigma$ and whose homsets are given by $\Hom(c, c') := \Path(c, c')$.

Given a category presentation $C$, let $\semp{C}\coloneqq \Free(C_\Sigma)/C_E = \Free(C_\Sigma)/\approx_C$. We call $\semp{C}$ the \textbf{category presented by} $C$. If $p : c \to c'$ is a $C$-path, then we let $[p] \in \semp{C}(c,d)$ denote its corresponding $\approx_C$-equivalence class in $\semp{C}$.
\end{Def}

Although we are interested in general category presentations (together with instance and profunctor presentations on them, as we will see in the subsequent sections), for examples and counterexamples we will concentrate on the particular case in which category presentations have only one sort $*$, in which case we call them \textbf{monoid presentations}. Much of the subtlety of the theory is already present in this class of examples, and it frees us from having to keep track of sorts.

\begin{Ex}\label{example monoid presentation}
Given a monoid presentation $M$, we write $M = \langle \Fun(M) \mid M_E \rangle$ following the usual notation for e.g. group presentations. For instance, let $M = \langle f, g \mid f.g = g.f\rangle$. Then $\semp{M}$ is the free commutative monoid on two generators. Moreover, let $F: M \to M$ be defined by $F(f) :\equiv f.f$, $F(g) :\equiv f.g$. This is a well defined morphism of category presentations because $F(f.g) \equiv f.f.f.g \approx_M f.g.f.f \equiv F(g.f)$.

We use the notation $\langle \Fun(M) \mid \varnothing \rangle$ to denote monoid presentations with no equations, and $\langle \varnothing \mid \varnothing \rangle$ to denote trivial monoid presentations, with no function symbols or equations. Notice that since this notation omits the unique sort, it technically only defines a presentation up to unique isomorphism. We usually denote the unique sorts as $*$ unless specified otherwise.
\end{Ex}

We have made use of congruences of categories in order to define the semantics of category presentations. We will now use them a second time in order to define a quotient of the ``bigger'' category $\ncat{CatPr}$.

\begin{Lemmarep}\label{lem presentation morphism preserves approx}
Given a morphism $F: C \to C'$ of category presentations, if $p_0 \approx_C p_1$, then $F(p_0) \approx_{C'} F(p_1)$.
\end{Lemmarep}
\begin{proof}
Straightforward induction, see Rem. \ref{rem structural induction on inference rules}.
\end{proof}

\begin{Def}
We say that two morphisms $F, G:C\to D$ of category presentations are \textbf{provably equal}, and write $F \approx G$, if and only if for each $c \in \Sort(C)$, the sorts $F(c)$ and $G(c)$ are equal and for each $f \in \Fun(C)$, $F(f) \approx_D G(f)$. It is easily shown using Lemma \ref{lem presentation morphism preserves approx} that the relation $\approx$ is a congruence on $\ncat{CatPr}$ , so we define $\ncat{CatPr}_\approx\coloneqq\ncat{CatPr}/{\approx}$.
\end{Def}

The construction $\semp{-}$ from Definition \ref{def semantic category} extends to a functor $\semp{-}: \ncat{CatPr} \to \ncat{Cat}$. This functor sends provably equal morphisms to equal morphisms, so it induces a unique functor $\widetilde{\semp{-}}:\ncat{CatPr}_\approx\to\ncat{Cat}$ by Lemma \ref{lem quotient category universal property}. In what follows we will abuse notation and allow $\semp{-}$ to also denote $\widetilde{\semp{-}}$.

\begin{Prop} \label{prop cat presentations and categories are equiv}
The functor $\semp{-}:\ncat{CatPr}_\approx\to\ncat{Cat}$ described above is an equivalence of categories.
\end{Prop}

\begin{proofsketch}
Given a small category $\cat{C}$, let $\can(\cat{C})$ denote its canonical category presentation\footnote{Also see \cite[Remark 4.21]{schultz2017algebraic} for a similar construction for algebraic theories.}, whose sorts are the objects of $\cat{C}$, whose function symbols are the morphisms of $\cat{C}$, and whose set of equations is the set of pairs of paths in $\can(\cat{C})$ whose compositions in $\cat{C}$ are equal. This construction can be shown to extend to a functor $\can: \ncat{Cat} \to \ncat{CatPr}_\approx$ which is a quasi-inverse for $\semp{-}$.
\end{proofsketch}

Note that $\ncat{CatPr}$ is not equivalent to $\ncat{Cat}$, as the following example shows. On the contrary, equality of morphisms in $\ncat{CatPr}$ is a strict, purely syntactic notion.

\begin{Ex}\label{ex mistake in alg databases}
Let $N \coloneqq \langle s \mid \varnothing \rangle$ and $M \coloneqq \langle f, g \mid f.g = g.f \rangle$ denote monoid presentations. Let $F, G : N \to M$ denote morphisms of category presentations defined by $F(s) :\equiv f.g$ and $G(s) :\equiv g.f$. Now $\semp{F}([s]) = [f.g] = [g.f] =\semp{G}([s])$, so $\semp{F} = \semp{G}$, but $F\neq G$ as morphisms in $\ncat{CatPr}$, so $\semp{-}$ is not faithful.
\end{Ex}

\section{Profunctor Presentations} \label{section profunctor presentations}

In this section we introduce the notions of uncurried and curried profunctor presentations. First we recall the definition of profunctors and profunctor composition.

\begin{Def} \label{def profunctor}
Let $\ncat{2}$ denote the category with two objects $0$ and $1$, and a single non-identity morphism $\leq : 0 \to 1$. Given categories $\cat{C}$ and $\cat{D}$, a \textbf{profunctor} $\cat{P}:\cat{C}\slashedrightarrow\cat{D}$ is a category $\cat{P}$ equipped with a functor $\pi:\cat{P}\to\ncat{2}$, such that $\pi^{-1}(0)=\cat{C}$ and $\pi^{-1}(1)=\cat{D}$. We say a morphism $f$ in $\cat{P}$ is a \textbf{cross-morphism} if $\pi(f) = \leq$.
Given profunctors $\cat{P}, \cat{P}':\cat{C}\slashedrightarrow\cat{D}$ a morphism of profunctors is a functor $F: \cat{P}\to\cat{P'}$ which restricts to the identity on $\cat{C}$ and $\cat{D}$.%
\footnote{There is a more general notion of morphism between possibly non-parallel profunctors, obtained from identifying profunctors with the category $\ncat{Cat}/\ncat{2}$ of small categories over $\ncat{2}$. We postpone discussion of this notion of morphism to Section \ref{section double categories}.} Let $\ncat{Prof}(\cat{C}, \cat{D})$ denote the resulting category.
\end{Def}

\begin{Rem} \label{rem profunctors as functors}
A profunctor $\cat{P}: \cat{C} \slashedrightarrow \cat{D}$ can equivalently be defined as a functor $\cat{P}: \cat{C}^{\op} \times \cat{D} \to \ncat{Set}$, see \cite[Example 2.19]{schultz2017algebraic} and \cite{Joyal2022}. Moreover, the Cartesian closure of $\ncat{Cat}$ provides an equivalence $\ncat{Prof}(\cat{C}, \cat{D}) \simeq [\cat{C}^\op \times \cat{D}, \ncat{Set}] \simeq [\cat{C}^\op, \ncat{Set}^\cat{D}]$ which allows us to think of a profunctor $\cat{P}$ also as a functor $\cat{P}: \cat{C}^\op \to \ncat{Set}^\cat{D}$. Each of these ways of thinking about profunctors suggests a different syntactic presentation.
\end{Rem}

\begin{Def} \label{def profunctor composition}
Given profunctors $\cat{P} : \cat{C} \xslashedrightarrow{} \cat{D}$ and $\cat{Q}: \cat{D} \xslashedrightarrow{} \cat{E}$, let $\cat{P} \odot \cat{Q}: \cat{C} \xslashedrightarrow{} \cat{E}$ denote the profunctor defined objectwise by the coend
\begin{equation}
    (\cat{P} \odot \cat{Q})(c,e) = \int^{d \in \cat{D}} \cat{P}(c,d) \times \cat{Q}(d, e).
\end{equation}
We call this the \textbf{composite profunctor} of $\cat{P}$ and $\cat{Q}$.
\end{Def}

We now give a concrete description of the sets $(\cat{P} \odot \cat{Q})(c,e)$, for $c \in \cat{C}$ and $e \in \cat{E}$, which we will take as our concrete definition of $(\cat{P} \odot \cat{Q})$. Let $\simeq_0$ denote the relation on $\sum_{d \in \cat{D}} \cat{P}(c,d) \times \cat{Q}(d,e)$ defined by $(p, \cat{Q}(g,1)(q')) \simeq_0 (\cat{P}(1,g)(p), q')$ for all morphisms $g : d \to d'$ in $\cat{D}$, $p \in \cat{P}(c,d)$, and $q' \in \cat{Q}(d',e)$.

Let $\simeq$ denote the smallest equivalence relation containing $\simeq_0$.
Then we set
\begin{equation} \label{eq description of profunctor composition hom}
    (\cat{P} \odot \cat{Q})(c,e) = \left( \sum_{d \in \cat{D}} \cat{P}(c,d) \times \cat{Q}(d,e) \right)/{\simeq}.
\end{equation}

We let $\langle p,q \rangle$ denote the equivalence class of the pair $(p,q)$ with $p \in \cat{P}(c,d)$ and $q \in \cat{Q}(d,e)$. If $f : c' \to c$ and $h: e \to e'$ are morphisms in $\cat{C}$ and $\cat{E}$ respectively, then
$(\cat{P} \odot \cat{Q})(f,h)(\langle p,q \rangle) = \langle \cat{P}(f, 1)(p),\cat{Q}(1, h)(q) \rangle$. 

It is not hard to show that $\odot$ defines a weakly associative and weakly unital\footnote{The unit profunctor on a category $\cat{C}$ is given by the Hom functor $\text{Hom}_{\cat{C}} : \cat{C}^\op \times \cat{C} \to \ncat{Set}$.} composition operation. This allows us to think of profunctors as the horizontal arrows of a double category. We will explore this idea more fully in Section \ref{section double categories}.

\subsection{Uncurried Profunctor Presentations}

\begin{Def}\label{def uncurried profunctor presentation}
Given category signatures $\Sigma_0$ and $\Sigma_1$, an \textbf{uncurried profunctor signature} $\Pi$ from $\Sigma_0$ to $\Sigma_1$ consists of a set $\Fun(\Pi)$ whose elements we call \textbf{profunctor function symbols}, along with two functions $s:\Fun(\Pi)\to\Sort(\Sigma_0)$, $t:\Fun(\Pi)\to\Sort(\Sigma_1)$. For $x\in\Fun(\Pi)$, we write $x : c \to d$ to mean that $s(x)=c$ and $t(x)=d$. We define the \textbf{category signature} $|\Pi|$ \textbf{associated to} $\Pi$ by setting $\Sort(|\Pi|)\coloneqq \Sort(\Sigma_0)+\Sort(\Sigma_1),$ and $\Fun(|\Pi|)\coloneqq \Fun(\Sigma_0)+\Fun(\Pi)+\Fun(\Sigma_1).$

Given a path $p: c \to d$ in $|\Pi|$, we say $p$ is a \textbf{cross-path} if $c \in \Sort(\Sigma_0)$ and $d \in \Sort(\Sigma_1)$. Let $\CPath(\Pi)$ denote the set of cross-paths of $\Pi$.
By an \textbf{uncurried profunctor equation} on $\Pi$, we mean an equation $p=q : c \to d$ between cross-paths $p: c \to d$ and $q: c\to d$ of $\Pi$. Given category presentations $C$ and $D$, a $(C,D)$-\textbf{uncurried profunctor presentation} $P$ is a pair $P = (P_\Pi, P_E)$ consisting of an uncurried profunctor signature $P_\Pi$ from $C_\Sigma$ to $D_\Sigma$ and a set $P_E$ of uncurried profunctor equations over $P_\Pi$. Let $\Fun(P)\coloneqq \Fun(P_\Pi)$, and $\CPath(P)\coloneqq \CPath(P_\Pi)$. We define the \textbf{associated category presentation} $|P|$ by $|P|_\Sigma \coloneqq |P_\Pi|$ and $|P|_E\coloneqq C_E + P_E + D_E$.  We write $p=_P q$ to mean that $p=q$ is an equation in $P_E$.  We define the relation $\approx_P$ as the restriction of $\approx_{|P|}$ to $\CPath(P)$.  (See Lemma \ref{lem uncurried approx characterization} for an inductive characterization of $\approx_P$.)
\end{Def}



\begin{Def}\label{def morphism of uncurried profunctor presentations}
Given uncurried profunctor presentations $P, P': C \pro D$, a \textbf{morphism of uncurried profunctor presentations} $F : P \to P'$ is 
a function $F:\Fun(P_\Pi)\to\CPath(P'_\Pi)$ such that if $p: c \to d$, then $F(p): c \to d$, and such that $|F|:|P|\to |P'|$ is a morphism of category presentations, where $|F|:|P|\to |P'|$ is the extension of $F$ to $|P|$ acting as the identity on $C$ and $D$.  We extend $F$ to cross-paths by $F(p):\equiv |F|(p)$ (using Definition \ref{def morphism of category presentations}).  We say that two such morphisms $F$ and $F'$ are \textbf{provably equal}, and write $F\approx F'$, if $|F|\approx |F'|$.

Let $\ncat{UnCurr}(C,D)$ denote the category of uncurried profunctor presentations from $C$ to $D$, where the composition $G \circ F$ of morphisms $F: P \to P'$ and $G: P' \to P''$ is defined as the function $\Fun(P_\Pi) \to \CPath(P''_\Pi)$ given by $p \mapsto |G|(F(p))$.
Let $\ncat{FinUnCurr}(C,D)$ denote the corresponding full subcategory of finite uncurried profunctor presentations, and let $\ncat{UnCurr}_{\approx}(C,D) \coloneqq \ncat{UnCurr}(C,D)/{\approx}$.
\end{Def}

\begin{Def}\label{def semantics of uncurried profunctor presentations}
Given an uncurried profunctor presentation $P : C \pro D$, we let $\semp{P}$ denote the profunctor $\semp{|P|}$, with $\pi:\semp{|P|}\to\ncat{2}$ sending all sorts of $C$ to $0$ and all sorts of $D$ to $1$. This construction extends to a functor $\semp{-}:\ncat{UnCurr}(C,D)\to\ncat{Prof}(C,D)$ for all $C$ and $D$.
\end{Def}

The following result can be proved using a similar argument as in the proof of Proposition \ref{prop cat presentations and categories are equiv}.

\begin{Prop}\label{prop equivalence between uncurr and prof}
For all category presentations $C$ and $D$, the induced functor $\semp{-}:\ncat{UnCurr}_{\approx}(C,D) \to \ncat{Prof}(\semp{C},\semp{D})$ is an equivalence.
\end{Prop}
\begin{proofsketch}
    As in the proof sketch for Proposition \ref{prop cat presentations and categories are equiv}, we define a quasi-inverse $\text{Can}: \ncat{Prof}(\semp{C},\semp{D}) \to \ncat{UnCurr}_{\approx}(C,D)$ by assigning to each profunctor $\cat{P}$ a canonical uncurried presentation between $C$ and $D$ whose set of profunctor symbols is the set of cross-morphisms of $\cat{P}$ and whose profunctor equations consist of pairs of paths in $\text{Can}(\cat{P})$ whose compositions in $\cat{P}$ are equal.
\end{proofsketch}

\begin{Def}
Given category presentations $C$ and $D$, we say that a profunctor $\cat{P}: \semp{C} \slashedrightarrow \semp{D}$ is \textbf{finitely uncurried presentable} if there exists a finite uncurried profunctor presentation $P : C \slashedrightarrow D$ and an isomorphism $\cat{P} \cong \semp{P}$ in $\ncat{Prof}(\semp{C}, \semp{D})$.
\end{Def}

Uncurried profunctor presentations are natural from the point of view of profunctors as categories over $\ncat{2}$. However, if we restrict to \emph{finite} presentations, then they are not suitable for presenting profunctor composition, as the following result illustrates.

\begin{Prop} \label{prop uncurried prof don't admit finite composition}
The class of finitely uncurried presentable profunctors is not closed under composition.
\end{Prop}

\begin{proof}
Consider monoid presentations $C = \langle \varnothing \mid \varnothing \rangle$, $D = \langle f \; \mid \varnothing \rangle $ and $E = \langle \varnothing \mid \varnothing \rangle$. For clarity, we refer to their unique sorts by $c, d$ and $e$, respectively.
Let $P: C \pro D$ be the uncurried profunctor presentation with one profunctor function symbol
$p: c \to d$
and no equations, and
let $Q : D \pro E$ be the uncurried profunctor presentation with one profunctor function symbol $q : d \to e$ and no equations.
Then $(\semp{P} \odot \semp{Q})(c,e) = \{ [p.f^n.q] \mid n \geq 0 \}$ is an infinite set.
Now suppose that $R: C \slashedrightarrow E$ is an uncurried profunctor presentation such that $\semp{R} \cong \semp{P} \odot \semp{Q}$. Then in particular $\semp{R}(c,e) \cong (\semp{P} \odot \semp{Q})(c,e)$ is infinite. But $C$ and $E$ have no function symbols which implies that $R$ must have infinitely many profunctor function symbols.

\end{proof}

Proposition \ref{prop uncurried prof don't admit finite composition} poses a potential problem for applications of computational categorical algebra. For instance, this is a problem in categorical database theory, where profunctor presentations are used to represent database queries \cite{schultz2017integration}. Curried profunctor presentations, which we introduce later in this section, do not share this deficiency. In order to define them, though, we first introduce instance presentations.

\subsection{Instance Presentations}

\begin{Not}
    Given a category $\cat{C}$, we will refer to a functor $\inst{I} : \cat{C} \to \ncat{Set}$ as an \textbf{instance} on $\cat{C}$, or $\cat{C}$-\textbf{instance}. We also write $\cat{C}\hyp\ncat{Inst} := \ncat{Set}^\cat{C}$.
\end{Not}



An instance $\cat{C} \to \ncat{Set}$ is equivalently a functor $\inst{I} : \ncat{1}^\op\times\cat{C}\to\ncat{Set}$, i.e. a profunctor $\inst{I} : \ncat{1}\slashedrightarrow\cat{C}$. This allows us to define instance presentations as a particular case of uncurried profunctor presentations.

\begin{Def}\label{def instance presentation}
Let $1$ be the category presentation with a single sort $*$ and no function symbols or equations. Given a category presentation $C$, a \textbf{$C$-instance signature} is an uncurried profunctor signature from $1$ to $C$, and a \textbf{$C$-instance presentation} is an uncurried profunctor presentation from $1$ to $C$. Let $C\hyp\ncat{InstPr} \coloneqq \ncat{UnCurr}(1,C)$ and $C\hyp\ncat{InstPr}_\approx\coloneqq \ncat{UnCurr}_\approx(1,C)$.
\end{Def}

Notice that the semantics of a $C$-instance presentation $I$ is a functor $\semp{I}: \semp{1}^\op \times \semp{C} \to \ncat{Set}$. Since $\semp{1}^\op \cong \ncat{1}$, it is equivalently a $\semp{C}$-instance $\semp{I} : \semp{C} \to \ncat{Set}$.

\begin{Not}
Given a $C$-instance presentation $I$, we call the profunctor function symbols of $I$ \textbf{generators}, writing $\text{Gen}(I) \coloneqq \Fun(I)$, and we call cross-paths of $I$ \textbf{terms}, writing $\Term(I)\coloneqq \CPath(I)$. Note that $|I|_E = I_E + C_E$. 
Since there is only one sort in $1$, we abbreviate $(x : * \to c) \in \Gen(I)$ to $(x : c) \in \Gen(I)$, to be read as ``$x$ is a term in $I$ of type $c$''. As in Definition \ref{def morphism of uncurried profunctor presentations}, $F$ also acts on terms by $F(t) :\equiv |F|(t)$.
\end{Not}

\begin{Ex}\label{example monoid action presentation}
Given a monoid presentation $M$, we refer to a $M$-instance presentation $I$ as a \textbf{$M$-action presentation}, and we write $I = \langle \Gen(I) \mid I_E \rangle_M$.
For example, let $N = \langle s \mid \varnothing \rangle$ be a presentation for the free monoid on one generator, i.e $\semp{N} \simeq \mathbb{N}$. Then $I = \langle x,y \mid x.s = x \rangle_N$ presents the set $\semp{I}(*) = \{[x], [y], [y.s], [y.s.s], \dots\}$ together with the $\semp{N}$-action $\semp{I}(s)([y.s^k]) = [y.s^k] \cdot s \coloneqq [y.s^{k+1}]$ and $[x] \cdot s \coloneqq [x]$. This can be identified with the set $\mathbb{N} \sqcup \{\infty\}$ and the $\mathbb{N}$-action $n \cdot s \coloneqq n+1$ and $\infty \cdot s \coloneqq \infty$.
\end{Ex}


\begin{toappendix}
For the subsequent proofs of this section and the next, we must prove many statements of the form $s \approx_P t \Rightarrow A(s,t)$ for terms $s$ and $t$ in a $(C,D)$-uncurried profunctor presentation $P$, and where $A(s,t)$ is a statement involving $s$ and $t$. We could attempt to do this by proving $s\approx_{|P|} t \Rightarrow A(s,t)$ by induction, but this would require us to extend $A(s,t)$ to arbitrary paths $s$ and $t$ which are not necessarily cross-paths.  The following result provides a useful inductive characterization for $\approx_P$ in its own right.


\begin{Lemma}\label{lem uncurried approx characterization}
Given a $(C,D)$-uncurried profunctor presentation $P$, the relation $\approx_P$ coincides with the relation $\sim$ on cross-paths of $P$ defined inductively by the following inference rules:
\begin{enumerate}
        \item For any equation $s =_P t$, $s \sim t$.
        \item For any cross-path $t$ of $P$, $t \sim t$.
        \item If $s \sim t$, then $t \sim s$.
        \item If $r\sim s$ and $s\sim t$, then $r\sim t$.
        \item For any $(g:d\to d')\in\Fun(D)$ and cross-paths $s, t : c\to d$, if $s \sim t$, then $s.g\sim t.g$.
        \item For any $(f:c'\to c)\in\Fun(C)$ and cross-paths $s, t : c\to d$, if $s \sim t$, then $f.s\sim f.t$.
        \item If $t:c\to d$ is a cross-path of $P$ and $g=_D g' : d\to d'$, then $t.g\sim t.g'$.
        \item If $t:c\to d$ is a cross-path of $P$ and $f=_C f' : c'\to c$, then $f.t\sim f'.t$.
\end{enumerate}
\end{Lemma}

\begin{proof}
A straightforward induction on the structure of $\sim$ (analogous to Remark \ref{rem structural induction on inference rules}) shows that $\sim$ is contained in $\approx_{|P|}$, and hence in $\approx_{P}$. To show the converse, let us show that if $s\approx_{|P|}t$, then either $s$ and $t$ are $C$-paths and $s\approx_C t$, $s$ and $t$ are $D$-paths and $s\approx_D t$, or $s$ and $t$ are cross-paths of $I$ and $s\sim t$. Let $A(s,t) = $ ``either $s$ and $t$ are $C$-paths and $s\approx_C t$, $s$ and $t$ are $D$-paths and $s\approx_D t$, or $s$ and $t$ are cross-paths of $I$ and $s\sim t$''. Now we will prove that $s \approx_{|P|} t$ implies $A(s,t)$ by induction on $\approx_{|P|}$, given by Definition \ref{def approx characterization}. If $s =_{|P|} t$, then either $s =_C t$ if they are $C$-paths, $s =_D t$ if they are $D$-paths, or $s =_P t$ if they are cross-paths, thus $A(s,t)$ holds.  Since $A$ is an equivalence relation, rules 2-4 preserve $A$ (see Remark \ref{rem structural induction on inference rules}).

For rule 5 of Definition \ref{def approx characterization}, suppose that $f:p'\to p''$ and $s,t :p\to p'$ are $|P|$-paths such that $s \approx_{|P|} t$ and $A(s,t)$.  If $p''$ is a sort of $C$, then clearly $s.f \approx_C t.f$ and $A(s.f,t.f)$ holds.  If $p$ is a sort of $D$, then clearly $s.f \approx_D t.f$ and $A(s.f, t.f)$ holds.  If $p$ is a sort of $C$ and $p'$ and $p''$ are sorts of $D$, then $s$ and $t$ are cross-paths and $s\sim t$.  Thus it follows that $s.f\sim t.f$ by rule 5 above, hence $A(s.f, t.f)$.  Finally, we consider the case where $p$ and $p'$ are sorts of $C$ and $p''$ is a sort of $D$.  In this case we have $s\approx_C t$ and $f$ is a cross-path.  We show that $s.f\sim t.f$, and hence $A(s.f, t.f)$, by a nested induction on $s\approx_C t$.  More specifically, we show that whenever $s\approx_C t : c'\to c$, $B(s,t)$ holds, where $B(s,t)\coloneqq \text{``for all cross-paths $f:c\to d$ we have $s.f\sim t.f$''}$.  If $s=_C t$, then $s.f\sim t.f$ follows for all cross-paths $f$ by rule 8 above.  Clearly $B$ is an equivalence relation, so the verifications for rules 2-4 of Def. \ref{def approx characterization} are automatic.  For rule 5 of Def. \ref{def approx characterization}, we assume that $B(s,t)$ holds for $s,t:c''\to c'$ with $s\approx_C t$ and we let $(u:c'\to c)\in\Fun(C)$.  Then for any cross-path $f:c\to d$, we obtain $s.u.f\sim t.u.f$ by instantiating $B$ on the cross-path $u.f$.  For rule 6 of Def. \ref{def approx characterization}, we assume that $B(s,t)$ holds for $s,t:c'\to c$ with $s\approx_C t$ and we let $(u:c''\to c')\in\Fun(C)$.  Then for any cross-path $f:c\to d$ we obtain $u.s.f\sim u.t.f$ by applying rule 6 above to the instantiation of $B(s,t)$ at $f$.

The verification of rule 6 of Def. \ref{def approx characterization} is dual to the verification of rule 5 that we have just completed.
\end{proof}

The following lemma follows immediately as an important special case of Lemma \ref{lem uncurried approx characterization}.

\begin{Lemma}\label{lem instance approx characterization}
Given a $C$-instance presentation $I$, the relation $\approx_I$ coincides with the relation $\sim$ on terms of $I$ defined inductively by the following inference rules:
\begin{enumerate}
        \item For any equation $s =_I t$, $s \sim t$.
        \item For any term $t$ of $I$, $t \sim t$.
        \item If $s \sim t$, then $t \sim s$.
        \item If $r\sim s$ and $s\sim t$, then $r\sim t$.
        \item For any $C$-function symbol $f:c\to d$ and terms $(s, t : c)$, if $s \sim t$, then $s.f\sim t.f : d$.
        \item If $t:c$ is a term of $I$ and $(p=_C q : c\to d)$, then $t.p\sim t.q:d$.
\end{enumerate}
\end{Lemma}

\end{toappendix}

\begin{Rem}
A profunctor $\cat{P}: \cat{C} \slashedrightarrow \cat{D}$ can be equivalently be thought of as a functor $\cat{P}: \cat{C}^{\op} \times \cat{D} \to \ncat{Set}$. This could be thought to give rise to a different notion of presentation: indeed, it is possible to define an involution $C \mapsto C^\op$ and an operation $(C, D) \mapsto C \times D$ in $\ncat{CatPr}$ such that one may define a profunctor presentation as a $(C^\op \times D)$-instance presentation. However, $(C^\op \times D)$-instance presentations turn out to be equivalent to $(C,D)$-uncurried profunctor presentations.

In contrast, if we use the Cartesian closure of $\ncat{Cat}$ to think of a profunctor $\cat{P}: \cat{C} \pro \cat{D}$ as a functor $\cat{P}: \cat{C}^{\op} \to \cat{D}\hyp\ncat{Inst}$, it becomes intuitive to define a profunctor presentation $P: C \pro D$ as a diagram of $D$-instance presentations, i.e. a family of $D$-instance presentations and morphisms between them indexed by the sorts and function symbols of $C$. We will show that in this case the resulting notion turns out to have starkly different properties. This might be surprising since the equivalence of categories $\ncat{Prof}(\cat{C}, \cat{D}) \simeq [\cat{C}^{\op}, \cat{D}\hyp\ncat{Inst}]$ is widely used throughout the categorical literature, often without mention.
\end{Rem}

\subsection{Curried Profunctor Presentations} \label{section curried presentations}


\begin{Def}\label{def curried profunctor presentation}
Given category presentations $C$ and $D$, a \textbf{curried profunctor presentation} $P : C \xslashedrightarrow{} D$ consists of an assignment of a $D$-instance presentation $P(c)$ to every sort $c \in C$, and an assignment of a $D$-instance presentation morphism $P(f): P(c') \to P(c)$ to every function symbol $(f: c \to c') \in \Fun(C)$, such that for each equation $p =_C p'$ we have $P(p) \approx P(p')$, where $P(p)$ is defined as $P(p) \coloneq P(p_n) \circ \dots \circ P(p_1)$ for $p \equiv p_1.\cdots.p_n$. We say $P: C \pro D$ is \textbf{finite} if $C$ and $D$ are finite and $P(c)$ is finite for all $c$.
\end{Def}

\begin{Ex}\label{example curried profunctor presentation}
    Given monoid presentations $M$ and $N$, we refer to an (un)curried profunctor presentation $P: M \pro N$ as a \textbf{(un)curried $(M,N)$-bimodule presentation}.
    For example, consider the previously introduced monoid presentations $M = \langle f, g \mid f.g = g.f \rangle$ and $N = \langle s \mid \varnothing \rangle$. Then a curried $(M, N)$-bimodule presentation is given by a choice of an $N$-action presentation $P(*)$ together with endomorphisms $P(f), P(g): P(*) \to P(*)$ such that $P(g) \circ P(f) \approx P(f) \circ P(g)$. Let $P(*)$ be the $N$-action presentation of Example \ref{example monoid action presentation}, $P(*) := \langle x, y \mid x.s = x\rangle_N$. As for the endomorphisms, not only must they commute (up to provable equality), but they must also respect the fact that $x$ is a fixed point (up to provable equality) of the action. This means that they must both send $x$ to a term provably equal to $x$. The reader may check that a valid choice is as follows: $P(f)(x) :\equiv P(g)(x) :\equiv x$, $P(f)(y) :\equiv y.s$ and $P(g)(y) :\equiv y.s.s$.
\end{Ex}

Since curried profunctor presentations are essentially diagrams of instance presentations, morphisms between them can be thought of as presentations of natural transformations.

\begin{Def}\label{def morphism of curried presentations}
Given curried profunctor presentations $P, P' : C \pro D$, a \textbf{morphism of curried profunctor presentations} $F: P \to P'$ consists of a collection of $D$-instance presentation morphisms $F_{c}:P(c)\to P'(c)$, indexed by the sorts $c$ of $C$, such that for each function symbol $f: c \to c'$ in $C$, the following diagram commutes up to provable equality:
\begin{equation} \label{eq diagram of morphism of curried presentations}
\begin{tikzcd}[ampersand replacement=\&]
	{P(c')} \& {P'(c')} \\
	{P(c)} \& {P'(c)}
	\arrow[""{name=0, anchor=center, inner sep=0}, "{F_{c'}}", from=1-1, to=1-2]
	\arrow["{P(f)}"', from=1-1, to=2-1]
	\arrow["{P'(f)}", from=1-2, to=2-2]
	\arrow[""{name=1, anchor=center, inner sep=0}, "{F_{c}}"', from=2-1, to=2-2]
	\arrow["\approx"{description}, draw=none, from=0, to=1]
\end{tikzcd}
\end{equation}
It follows that \eqref{eq diagram of morphism of curried presentations} still commutes when $f$ is a path, by induction on $f$; thus these morphisms can be composed in the natural way. Let $\ncat{Curr}(C,D)$ denote the category of curried profunctor presentations from $C$ to $D$, and if $C$ and $D$ are finite, let $\ncat{FinCurr}(C,D)$ denote the full subcategory of finite curried profunctor presentations.

Given $P, P': C \pro D$, we say that two curried profunctor presentation morphisms $F,G:P\to P'$ are \textbf{provably equal}, and write $F\approx G$, if $F_c\approx G_c$ for each $c\in C$.
\end{Def}

\begin{Def}\label{def semantics of curried profunctor presentations, globular case}

Given a curried profunctor presentation $P : C \xslashedrightarrow{} D$, let $\semp{P}$ denote the functor $\semp{P}: \semp{C}^\op \to \semp{D}\hyp\ncat{Inst}$ given by $\semp{P}(c) \coloneq \semp{P(c)}$ for each $c \in \semp{C}$ and $\semp{P}([f]) \coloneq \semp{P(f_n)} \circ \dots \circ \semp{P(f_1)}$ for each morphism $[f] = [f_1.\cdots.f_n]$ in $\semp{C}$. Notice that well definition follows from the definition of curried profunctor presentation. Moreover, $\semp{P}$ thus defined is equivalently a profunctor $\semp{P} : \semp{C}^\op \times \semp{D} \to \ncat{Set}$ with $\semp{P}(c,d) = \semp{P(c)}(d)$. It is easy to see that this extends to a functor $\semp{-}: \ncat{Curr}(C,D) \to \ncat{Prof}(\semp{C},\semp{D})$.

We say a profunctor $\cat{P}: \semp{C} \pro \semp{D}$ is \textbf{finitely curried presentable} if there exists a finite curried profunctor presentation $P: C \pro D$ such that $\semp{P} \cong \cat{P}$.
\end{Def}

\begin{Rem}\label{rem bicategories of uncurried profunctor presentations}
    A curried profunctor presentation $P: C \slashedrightarrow D$ does not determine a functor $\semp{C}^\op\to D\hyp\ncat{InstPr}$. That would require the diagram (\ref{eq diagram of morphism of curried presentations}) to commute strictly, rather than up to provable equality. Now $P$ does determine a functor $P : \semp{C}^\op \to D\hyp\ncat{InstPr}_\approx$, but the converse is not true, because the morphisms $P(f)$ are actual data, rather than equivalence classes of such data.

    Given a category $\cat{C}$ and a congruence $\sim$ on $\cat{C}$, we can also consider the pair $(\cat{C},\sim)$ as a bicategory where given morphisms $f,g$ in $\cat{C}$, there is a unique 2-morphism $f \Rightarrow g$ if and only if $f \sim g$. With this definition $P$ defines a pseudofunctor $\semp{C}^\op \to (D\hyp\ncat{InstPr},\approx)$. We leave the development of the consequences of this observation to future work.
\end{Rem}

It is not hard to show that every profunctor can be presented by an uncurried profunctor presentation or a curried profunctor presentation, using ideas similar to the proof of Proposition \ref{prop cat presentations and categories are equiv}. However, if we restrict to \textit{finite} uncurried and curried profunctor presentations, then the situation is very different, as the following result shows.

\begin{Lemma}\label{lem finite uncurried but not finite curried}
There exist profunctors that are finitely uncurried presentable but not finitely curried presentable.
\end{Lemma}

\begin{proof}
As in the proof of Prop. \ref{prop uncurried prof don't admit finite composition}, let $D = \langle f \mid \varnothing \rangle$ and $E = \langle \varnothing \mid \varnothing \rangle$ denote monoid presentations with unique sorts $d$ and $e$ respectively, and let $Q : D \slashedrightarrow E$ be the uncurried bimodule presentation given by $\Fun(Q) = \{ q : d \to e \}$ and no equations. Then $\semp{Q}(d,e) = \{[p], [f.p], [f.f.p], \dots \}$ is an infinite set, and $\semp{Q}(f)([f^n.p]) = [f^{n+1}.p]$. If $Q' : D \slashedrightarrow E$ is a curried bimodule presentation for $\semp{Q}$, then $Q'(d)$ must be a $E$-instance presentation such that $\semp{Q'(d)}(e) \cong \semp{Q}(d,e)$. This implies that $Q'(d)$ has infinitely many terms. Since $E$ has no function symbols, $Q'(d)$ must have infinitely many generators. Thus it cannot be a finite curried bimodule presentation.
\end{proof}

In contrast to Proposition \ref{prop uncurried prof don't admit finite composition}, we define a notion of syntactic composition for curried profunctor presentations that is semantically correct.

\begin{Def} \label{def curried composition instance signature}
Given category presentations $C$, $D$, and $E$, curried profunctor presentations $P : C \slashedrightarrow D$ and $Q : D \slashedrightarrow E$, and a sort $c \in C$, we define the $E$-instance presentation $(P\ccomp Q)(c)$ as follows. let $(P \ccomp Q)(c)_\Pi$ denote the $E$-instance signature whose set of generators consists of pairs of generators $(p : d) \in \Gen(P(c))$ and $(q : e) \in \Gen(Q(d))$, which we write as $(p \otimes q : e)$. We now extend the operation $\otimes$ to accept a pair $(s,t)$ of terms with $(s : d) \in \Term(P(c))$, $(t : e) \in \Term(Q(d))$ and return a term in $(P \ccomp Q)(c)$. Firstly, if $(p : d)\in \Gen(P(c))$ and $(t \equiv q.h : e) \in \Term(Q(d))$, where $q$ is a generator of $Q(d)$ and $h$ is a path in $E$, then set $p\otimes t :\equiv (p\otimes q).h$. Secondly, if $(s \equiv p.g : d)\in \Term(P(c))$ and $(t : e) \in \Term(Q(d))$, then set $s\otimes t :\equiv p\otimes Q(g)(t)$.

Let $(P\ccomp Q)(c)$ be the $E$-instance presentation with signature $(P\ccomp Q)(c)_\Pi$ and set of equations $(P\ccomp Q)(c)_E$ consisting of equations of the form
\begin{itemize}
    \item $(s\otimes q =_{(P\ccomp Q)(c)} s'\otimes q)$, where $(s =_{P(c)} s':d)$ and $(q:e)\in \Gen(P(d))$, or 
    \item $(p \otimes t =_{(P\ccomp Q)(c)} p\otimes t')$, where $(p:d) \in \Gen(P(c))$ and $(t =_{Q(d)} t':e)$.
\end{itemize}
\end{Def}

\begin{Lemmarep} \label{lem composition preserves provable equality of factors}
Given curried profunctor presentations $P : C \pro D$ and $Q : D \pro E$, the operation $\otimes$ respects provable equality. In other words, if $(s\approx_{P(c)}s':d)$ and $(t\approx_{Q(d)}t' :e)$, then $(s\otimes t \approx_{(P\ccomp Q)(c)} s'\otimes t' : e)$.
\end{Lemmarep}
\begin{proof}
First let us prove the following statements: 
\begin{enumerate}
    \item if $(t\approx_{Q(d)}t':e)$, then for all $(p:d)\in \Gen(P(c))$, $(p\otimes t \approx_{(P\ccomp Q)(c)} p\otimes t' :e)$.
    \item if $(t\approx_{Q(d)}t':e)$, then for all $(s :d)\in \Term(P(c))$, $(s \otimes t \approx_{(P\ccomp Q)(c)} s \otimes t' :e)$.
    \item if $(s\approx_{P(c)}s':d)$, then for all $(t :e)\in \Term(Q(d))$, $(s\otimes t\approx_{(P\ccomp Q)(c)}s'\otimes t:e)$.
\end{enumerate}
Given the definition of the equations in $P \ccomp Q$, it is straightforward to prove (1) by induction on $\approx_{Q(d)}$ using Lemma \ref{lem instance approx characterization}. Let us now consider (2). If $(s : d) \in \Term(P(c))$, then $s \equiv p.g$ for some generator $p$ and $D$-path $g: d' \to d$. Thus $s \otimes t \equiv (p \otimes Q(g)(t))$ and $s \otimes t' \equiv (p \otimes Q(g)(t'))$. Since $Q$ is a curried profunctor presentation and $t \approx_{Q(d)} t'$, we have $Q(g)(t) \approx_{Q(d')} Q(g)(t')$ by Lemma \ref{lem presentation morphism preserves approx}, and thus by (1), it follows that $s \otimes t \equiv (p \otimes Q(g)(t)) \approx_{(P \ccomp Q)(c)} p \otimes Q(g)(t') \equiv s \otimes t'$.

Now let us prove (3) by induction on $\approx_{P(c)}$ using Lemma \ref{lem instance approx characterization}.
If $(s =_{P(c)} s' : d)$ and $(t \equiv q.h: e) \in \Term(Q(d))$, then there is an equation $(s \otimes q =_{(P \ccomp Q)(c)} s'\otimes q: e)$, hence $s \otimes t \equiv (s\otimes q).h \approx_{(P \ccomp Q)(c)} (s'\otimes q).h \equiv s'\otimes t$ by rule 6 of Lemma \ref{lem instance approx characterization}. Rules 2-4 clearly follow from the fact that the property in the consequent of (3) is reflexive, symmetric and transitive.

For Rule 5, we want to show that if $s,s' \in \Term(P(c))$, $(s\approx_{P(c)}s':d)$ and for all $(t:e)\in \Term(Q(d))$, $(s\otimes t\approx_{(P\ccomp Q)(c)}s'\otimes t:e)$, then for any $D$-path $g: d \to d'$ and term $(t:e)\in \Term(Q(d'))$ we have $(s.g \otimes t \approx_{(P \ccomp Q)(c)} s'.g \otimes t : e)$. But this is equivalent to showing that $(s \otimes Q(g)(t) \approx_{(P \ccomp Q)(c)} s' \otimes Q(g)(t))$, which follows by the inductive hypothesis.

Finally, for Rule 6, suppose that $(g =_D g': d \to d')$. Then $Q(g) \approx Q(g')$ since $Q$ is a curried profunctor presentation. Thus for every $(s : d) \in \Term(P(c))$ and $(t: e) \in \Term(Q(d'))$, $s \otimes Q(g)(t) \approx_{(P \ccomp Q)(c)} s \otimes Q(g')(t)$ by (2) above, which is to say that $s.g \otimes t \approx_{(P\ccomp Q)(c)} s.g' \otimes t$.

Having proven statements (1), (2) and (3), the result follows by transitivity of $\approx_{(P \ccomp Q)(c)}$.


\end{proof}

\begin{Def}\label{def composite curried profunctor presentation}
Given curried profunctor presentations $P : C \slashedrightarrow D$ and $Q: D \slashedrightarrow E$ as above, we define the \textbf{composite presentation} $P \circledast Q$ as follows. We have already defined an $E$-instance $(P\circledast Q)(c)$ for each $c\in\Sort(C)$ in Definition \ref{def curried composition instance signature}. Given a function symbol $f: c \to c'$ in $C$, we define a function $(P \ccomp Q)(f) : \Gen((P \ccomp Q)(c')) \to \Term((P \ccomp Q)(c))$ by setting $(P \ccomp Q)(f)(p\otimes q):\equiv P(f)(p)\otimes q$.
\end{Def}


\begin{Lemmarep} \label{lem curried comp is well defined}
Given a $C$-function symbol $f: c \to c'$ in $C$, the function $(P \ccomp Q)(f)$ is a morphism of $E$-instances. Furthermore $P \ccomp Q$ is a well-defined curried profunctor presentation from $C$ to $E$.
\end{Lemmarep}

\begin{proof}
First notice that $(P\ccomp Q)(f)(s\otimes t) \equiv P(f)(s)\otimes t$ for all $(f:c\to c')\in\Path(C)$, $(s:d)\in\Term(P(c))$, and $(t:e)\in\Term(Q(d))$, as can be checked by a simple computation.

Now given a $C$-function symbol $f:c\to c'$, to show that $(P\ccomp Q)(f)$ is an $E$-instance morphism, we must show that it sends equations in $(P \ccomp Q)(c')$ to provable equalities in $(P \ccomp Q)(c)$. Since equations in $(P \ccomp Q)(c')$ have two possible forms, we consider them one by one. If $(s =_{P(c)} s':d)$ and $(q:e)\in\Gen(P(d))$, then $(P(f)(s)\otimes q \approx_{(P\ccomp Q)(c)} P(f)(s')\otimes q)$ by Lemma \ref{lem composition preserves provable equality of factors}. On the other hand, if $(p:d)\in\Gen(P(c))$ and $(t=_{Q(d)} t':d)$, then $(P(f)(p)\otimes t\approx_{(P\ccomp Q)(c)} P(f)(p)\otimes t')$ again by Lemma \ref{lem composition preserves provable equality of factors}.

Secondly, we must show that whenever $(f=_C f':c\to c')$ is an equation between $C$-paths, we have $(P \ccomp Q)(f) \approx (P\ccomp Q)(f')$. 
Thus this reduces to showing that $P(f)(p)\otimes q \approx_{(P \ccomp Q)(c)} P(f')(p)\otimes q$ for all generators $(p\otimes q)\in\Gen((P\ccomp Q)(c'))$. But this follows directly from the fact that $P(f)\approx P(f')$ and Lemma \ref{lem composition preserves provable equality of factors}.
\end{proof}

\begin{Lemmarep}\label{lem composition of curried presentations is functorial}
The assignment $(P, Q) \mapsto P \circledast Q$ defines a functor 
$$\ccomp : \ncat{Curr}(C,D) \times \ncat{Curr}(D,E) \to \ncat{Curr}(C,E).$$
\end{Lemmarep}
\begin{proof}
Suppose that $P, P': C \xslashedrightarrow{} D$ and $Q, Q' : D \xslashedrightarrow{} E$ are curried profunctor presentations and $\phi: P \to P'$ and $\psi: Q \to Q'$ are morphisms of curried profunctor presentations. Given $c \in C$, let $(\phi \ccomp \psi)_c : (P \ccomp Q)(c) \to (P' \ccomp Q')(c)$ denote the map of $E$-instance presentations defined by
$$(\phi \ccomp \psi)_c(p \otimes q) :\equiv (\phi_c(p) \otimes \psi_d(q))$$
for all $(p : d) \in \Gen(P(c))$ and $(q : e) \in \Gen(Q(d))$.
In what follows we make repeated use of Lemma  \ref{lem composition preserves provable equality of factors} without stating it explicitly.

First, notice that that
\begin{equation}\label{eq phi psi terms}
(\phi\ccomp\psi)_c(s\otimes t) \approx_{(P' \ccomp Q')(c)} \phi_c(s)\otimes \psi_d(t).
\end{equation}
for all $(s:d)\in \Term(P(c))$ and $(t:e)\in\Term(Q(d))$.
Indeed, expand $s\equiv p.g$ where $(p : d') \in \Gen(P(c))$ and $g:d' \to d$ is a $D$-path, and $Q(g)(t)\equiv q.h$. Then 
$$(\phi \ccomp \psi)_c(s \otimes t) \equiv (\phi \ccomp \psi)_c(p \otimes q.h) \equiv \phi_c(p) \otimes \psi_{d'}(q).h \approx_{(P'\ccomp Q')(c)} \phi_c(p)\otimes Q'(g)(\psi_d(t)) \equiv \phi_c(s) \otimes \psi_d(t)$$
where we have used ``naturality'' of $\psi$ as expressed by diagram \eqref{eq diagram of morphism of curried presentations}. The fact that $(\phi \ccomp \psi)_c$ is an $E$-instance presentation morphism then follows directly from \eqref{eq phi psi terms} by looking at the specific form of any equation in $P \ccomp Q$ and remembering that both $\phi_c$ and $\psi_d$ are morphisms of instance presentations.

Finally, let us show that the naturality condition expressed by diagram \eqref{eq diagram of morphism of curried presentations} holds for the family of morphisms $\{(\phi \ccomp \psi)_c\}_{c \in C}$. Given a function symbol $f:c\to c'$ in $C$ and a generator $(p\otimes q:e)\in \Gen((P\ccomp Q)(c'))$, where $(p:d) \in \Gen(P(c'))$,

\begin{align*}
    (\phi \ccomp \psi)_c \, ((P\ccomp Q)(f)(p \otimes q)) &\equiv \phi_c(P(f)(p)) \otimes \psi_d(q) \tag{def. of $\phi \ccomp \psi$}\\
    &\approx_{(P' \ccomp Q')(c)} P'(f)(\phi_{c'}(p)) \otimes \psi_d(q) \tag{naturality of $\phi$} \\
    &\equiv (P' \ccomp Q')(f) \, (\phi_{c'}(p) \otimes \psi_d(q)) \tag{def. of $(P' \ccomp Q')(f)$} \\
    &\equiv (P' \ccomp Q')(f)\, \left( (\phi \ccomp \psi)_{c'} (p \otimes q) \tag{def. of $\phi \ccomp \psi$} \right).
\end{align*}
%
\end{proof}

Now that we have seen that the composition operation $\circledast$ is well defined, we now show that its semantics are given by $\odot$. From the algorithmic point of view, this provides a proof of correctness of the composition algorithm given implicitly by $\circledast$.

\begin{Threp} \label{th curried composition semantics}
Given curried profunctor presentations $P: C \pro D$ and $Q: D \pro E$, there is an isomorphism
\begin{equation*}
    \mu^{P, Q} : \semp{P} \odot \semp{Q} \to \semp{P \ccomp Q}
\end{equation*}
of profunctors, natural in $P$ and $Q$. Furthermore, if $P : C \slashedrightarrow D$ and $Q: D \slashedrightarrow E$ are finite curried profunctor presentations, then $P \ccomp Q$ is a finite curried profunctor presentation. In particular, composition of profunctors preserves finite curried presentability.
\end{Threp}

\begin{proof}
It is obvious by construction that if $P : C \xslashedrightarrow{} D$ and $Q: D \xslashedrightarrow{} E$ are finite curried profunctor presentations, then $(P \ccomp Q)$ is a finite curried profunctor presentation. Thus all we need to do is show that there is an isomorphism of profunctors $\mu^{P,Q}: \semp{P} \odot \semp{Q} \cong \semp{P \ccomp Q}$ natural in $P$ and $Q$. We consider both profunctors as functors $\semp{C}^\op \times \semp{D} \to \ncat{Set}$. Hence, to give this isomorphism is to give a natural isomorphism, whose components are bijections $\mu^{P,Q}_{(c,e)} : (\semp{P} \odot \semp{Q})(c,e) \to \semp{P \ccomp Q}(c,e)$ for each $c \in C$ and $e \in E$. Moreover, we work with the concrete description of the sets $(\semp{P} \odot \semp{Q})(c,e)$ given by \eqref{eq description of profunctor composition hom}, hence these functions are defined on equivalence classes $\langle [s], [t]\rangle$ where $[s] \in \semp{P}(c,d)$ and $[t] \in \semp{Q}(d,e)$ for some $d \in D$. Notice that in this formulation of profunctors as functors, their action on morphisms are given by pre- and postcomposition, i.e. given $f: c' \to c$ a $C$-path and $k: e \to e'$ an $E$-path,
\begin{align*}
(\semp{P} \odot \semp{Q})([f],[k])(\langle[s], [t]\rangle) &= \langle \semp{P}([f])([s]), [t.k] \rangle \\
&= \langle [P(f)(s)], [t.k] \rangle \qquad \text{ and} \\
\semp{P \ccomp Q}([f], [k])([(p\otimes q).h]) &= [(P \ccomp Q)(f)(p \otimes q).h.k] \\
&= [((P)(f)(p) \otimes q).h.k].
\end{align*}

Let us then define, for each $c \in C$ and $e \in E$, the map $\mu^{P,Q}_{(c,e)}$ by the formula
$$\mu^{P,Q}_{(c,e)}(\langle [s], [t] \rangle) \coloneqq [s \otimes t]$$
for all $d \in D$, $(s: d) \in \Term(P(c))$ and $(t : e) \in \Term(Q(d))$. We will now prove that this is a well defined function. Notice that the domain of $\mu^{P,Q}_{(c,e)}$ is defined by two nested quotients by equivalence relations. The situation is clarified by the following diagram
\[\begin{tikzcd}[ampersand replacement=\&]
	{S_{(c,e)}} \& {\semp{P\ccomp Q}(c,e)} \& {(s,t)} \& {[s\otimes t]} \\
	{\widetilde{S}_{(c,e)}} \& {(\semp{P} \odot \semp{Q})(c,e)} \& {([s],[t])} \& {\langle[s],[t]\rangle}
	\arrow["F", from=1-1, to=1-2]
	\arrow["{q_1}"', from=1-1, to=2-1]
	\arrow["{q_2}"', from=2-1, to=2-2]
	\arrow["{\widetilde{F}}", dashed, from=2-1, to=1-2]
	\arrow["{\mu_{(c,e)}^{P,Q}}"', dashed, from=2-2, to=1-2]
	\arrow["{q_1}"', maps to, from=1-3, to=2-3]
	\arrow["{q_2}"', maps to, from=2-3, to=2-4]
	\arrow["F", maps to, from=1-3, to=1-4]
	\arrow[dashed, maps to, from=2-4, to=1-4]
\end{tikzcd}\]
where $S_{(c,e)} \coloneq \sum_{(s:d) \in \Term(P(c))} \{(t:e) \in \Term(Q(d))\}$ and $\widetilde{S}_{(c,e)} \coloneq \sum_{d \in D} \semp{P}(c,d)\times\semp{Q}(d,e)$, and the functions $F$, $q_1$ and $q_2$ are given by the the mappings on the right. We will first prove that $F$ is constant on the fibers of $q_1$, hence inducing a function $\widetilde{F}$ making the upper triangle in the diagram commute, and then that $\widetilde{F}$ is constant on the fibers of $q_2$, hence inducing our desired function $\mu_{(c,e)}^{P,Q}$. Both are relatively straightforward.

The first step amounts to showing that the expression $[s \otimes t]$ is a well defined function of the pair $([s], [t])$. Let $(s:d), (s' : d') \in \Term(P(c))$, $(t: e) \in \Term(Q(d))$ and $(t': e) \in \Term(Q(d'))$ be such that $([s], [t]) = ([s'], [t'])$. Then it must be the case that $d' = d$, $s \approx_{P(c)} s'$ and $t' \approx_{Q(d)} t$, from which $s \otimes t \approx_{(P\ccomp Q)(c)} s'\otimes t'$ by Lemma \ref{lem composition preserves provable equality of factors}.


For the second step, let $R \subseteq \widetilde{S}_{(c,e)} \times \widetilde{S}_{(c,e)}$ be the binary relation given by $([s],[t]) \, R \, ([s'],[t'])$ if and only if $\widetilde{F}([s],[t]) = [s\otimes t] = [s'\otimes t'] = \widetilde{F}([s'],[t'])$. Since the relation ``$q_2([s],[t])=q_2([s'],[t'])$'' is precisely the relation $\simeq$ from \eqref{eq description of profunctor composition hom}, and since $\simeq$ is defined as the smallest equivalence relation containing $\simeq_0$, it is enough to show that $R$
is an equivalence relation containing $\simeq_0$.
Since $R$ is clearly an equivalence relation, it only remains to show that $R$ contains $\simeq_0$.
Given $[s]: c \to d$, $[s']:c\to d'$, $[t]:d \to e$ and $[t']:d'\to e$, suppose $([s], [t]) \simeq_0 ([s'], [t'])$. Then there exists a $D$-path $f : d \to d'$ such that $[s.f] =  [s']$ and $[Q(f)(t')] = [t]$. Since $s.f \otimes t' \equiv s\otimes Q(f)(t')$, we have $\widetilde{F}([s],[t])=\widetilde{F}([s],[Q(f)(t')])=[s\otimes Q(f)(t')]=[s.f\otimes t']=\widetilde{F}([s.f],[t'])=\widetilde{F}([s'],[t'])$.

Thus we conclude that $\mu^{P,Q}_{(c,e)}$ is well defined. Meanwhile, naturality of $\mu^{P,Q}$ amounts to the fact that, for all $C$-paths $f: c' \to c$ and $E$-paths $k: e \to e'$, we have the following:
\begin{equation*}
\begin{tikzcd}[ampersand replacement=\&]
	{\langle[s],[t]\rangle} \&\& {\langle [P(f)(s)], [t.k]\rangle} \\
	{[s\otimes t]} \& {[(P\ccomp Q)(f)(s\otimes t).k]} \& {[P(f)(s) \otimes t.k]}
	\arrow["{\mu_{(c,e)}^{P,Q}}"', maps to, from=1-1, to=2-1]
	\arrow[maps to, from=1-1, to=1-3]
	\arrow["{\mu_{(c',e')}^{P,Q}}", maps to, from=1-3, to=2-3]
	\arrow[maps to, from=2-1, to=2-2]
	\arrow[Rightarrow, no head, from=2-2, to=2-3].
\end{tikzcd}
\end{equation*}

We have proved that $\mu^{P,Q}$ is a natural transformation. Let us now prove that it is an isomorphism. Fixing again $c$ and $e$, and setting $\alpha \coloneqq \mu^{P,Q}_{(c,e)}$ for the remainder of this proof, we define a map backwards $\beta : \semp{P \ccomp Q}(c,e) \to (\semp{P} \odot \semp{Q})(c,e)$ as follows.
If $x \in \semp{P \ccomp Q}(c,e)$, then it is of the form $x = [r]$ for some term $r \equiv (p \otimes q).h$ where
$(p : d) \in \Gen(P(c))$ and $(q.h : e) \in \Term(Q(d))$.
We then define $\beta(x) \coloneq \langle [p], [q.h] \rangle$. Let us show that $\beta$ is well defined by induction on $\approx_{(P \ccomp Q)(c)}$. 

For the base case, consider an equation $r =_{(P \ccomp Q)(c)} r'$. There are two possibilities:
\begin{itemize}
    \item $r \equiv s \otimes q$, $r' \equiv s' \otimes q$ where $s =_{P(c)} s': d$ and $(q:e) \in \Gen(Q(d))$. In this case, let $s \equiv p.g$ and $s' \equiv p'.g'$, hence $r \equiv p \otimes Q(g)(q)$ and $r' \equiv p' \otimes Q(g')(q)$. Then
    $$\beta([r]) = \langle [p], [Q(g)(q)]\rangle = \langle [p.g], [q] \rangle = \langle [p'.g'], [q] \rangle = \langle [p'], [Q(g')(q)]\rangle = \beta([r']).$$
    \item $r \equiv p \otimes t$, $r' \equiv p \otimes t'$ where $(p:d) \in \Gen(P(c))$ and $t =_{Q(d)} t':e$. In this case we conclude immediately that $\beta([r]) = \langle [p], [t] \rangle = \langle [p], [t'] \rangle = \beta([r'])$.
\end{itemize}
Thus the base case is proved. Meanwhile, rules 2-4 of Lemma \ref{lem instance approx characterization} are immediate from the fact that the relation $\{((p \otimes q).h, (p' \otimes q').h') \mid \langle [p], [q.h]\rangle = \langle [p'], [q'.h']\rangle\}$ is an equivalence relation.
For rule 5, suppose that $(r \approx_{(P \ccomp Q)(c)} r' : e)$, $k : e \to e'$ is an $E$-path and $\beta([r]) = \beta([r'])$.
Let $r \equiv (p\otimes q).h$ and $r' \equiv (p' \otimes q').h'$ so that $\langle [p], [q.h]\rangle = \langle [p'], [q'.h']\rangle$.
Then
$$\langle [p], [q.h.k]\rangle = (\semp{P}\odot\semp{Q})(1,[k])(\langle [p], [q.h]\rangle) = (\semp{P}\odot\semp{Q})(1,[k])(\langle [p'], [q'.h']\rangle) = \langle [p'], [q'.h'.k]\rangle.$$
Rule 6 follows by similar reasoning, from the fact that if $(k =_E k' : e \to e')$, then $[k] = [k']$ in $\semp{E}$ and hence $(\semp{P} \odot \semp{Q})(1, [k]) = (\semp{P} \odot \semp{Q})(1, [k'])$. Thus $\beta$ is well defined. 

Now it is easy to show that $\alpha$ and $\beta$ are inverse to each other. Indeed, $(\beta \circ \alpha)(\langle [p.g], [t] \rangle) = \beta([p.g \otimes t]) = \langle [p], [Q(g)(t)] \rangle = \langle [p.g], [t] \rangle$. Similarly, $(\alpha \circ \beta)([(p \otimes q).h]) = \alpha (\langle [p], [q.h] \rangle) = [p \otimes q.h] = [(p\otimes q).h]$.

Thus $\alpha = \mu^{P,Q}_{(c,e)}$ is a bijection for every $(c,e)$, and therefore $\mu^{P,Q}$ defines an isomorphism of profunctors.

Finally, we must show that $\mu : \semp{-} \odot \semp{=} \Rightarrow \semp{- \ccomp =}$ is a natural transformation between functors $\ncat{Curr}(C,D) \times \ncat{Curr}(D,E) \to \ncat{Prof}(\semp{C}, \semp{E})$. Let $\phi: P \to P'$ and $\psi: Q \to Q'$ be morphisms of curried profunctor presentations $P, P': C \pro D$ and $Q, Q': D \pro E$. Then naturality amounts to the verification
\begin{align*}
(\semp{\phi \ccomp \psi} \circ \mu^{P,Q}) (\langle [s], [t]\rangle) &= \semp{\phi \ccomp \psi} ([s \otimes t]) \\
&= [\phi_c(s) \otimes \psi_d(t)] \\
&= \mu^{P',Q'}(\langle [\phi_c(s)], [\psi_d(t)]\rangle) \\
&= (\mu^{P',Q'} \circ (\semp{\phi}\odot\semp{\psi})) (\langle[s],[t]\rangle).
\end{align*}
\end{proof}

Notice that we can extract an algorithm for computing the composition of two finite curried profunctor presentations by inspection of Definition \ref{def composite curried profunctor presentation}. In this sense, Theorem \ref{th curried composition semantics} proves the correctness of this algorithm with respect to the semantics.

\begin{Rem}
    Theorem \ref{th curried composition semantics} gives an alternative proof of Lemma \ref{lem finite uncurried but not finite curried}. Indeed, consider profunctors $\semp{P}$ and $\semp{Q}$ in Proposition \ref{prop uncurried prof don't admit finite composition}. It is not hard to see that $\semp{P}$ is finitely curried presentable and $\semp{P} \odot \semp{Q}$ is not finitely curried presentable. Hence, by Theorem \ref{th curried composition semantics} $\semp{Q}$ is not finitely curried presentable.
\end{Rem}


\begin{Ex} \label{ex curried example}
Here we give an example of curried profunctor presentation composition. Recall the $(M,N)$-bimodule presentation $P: M \pro N$ from Example \ref{example curried profunctor presentation}. Now let $O$ be the monoid presentation $O = \langle t \mid \varnothing \rangle$ and let $Q: N \pro O$ be the $(N,O)$-bimodule presentation $Q(*) = \langle q \mid q.t.t = q.t \rangle_O$ with $Q(s)(q) :\equiv q.t$.
Then the composite presentation $(P \ccomp Q) : M \pro O$ is given by
$$(P \ccomp Q)(*) = \langle x \otimes q, \ y \otimes q \mid (x \otimes q).t.t = (x \otimes q).t, \ (y \otimes q).t.t = (y \otimes q).t , \ (x\otimes q).t = x \otimes q\rangle_O$$
with 
\begin{itemize}
    \item $(P \ccomp Q)(f)(x \otimes q) \equiv P(f)(x) \otimes q \equiv x \otimes q$,
    \item $(P \ccomp Q)(f)(y \otimes q) \equiv P(f)(y) \otimes q \equiv y.s \otimes q \equiv y \otimes Q(s)(q) \equiv (y \otimes q).t$,
    \item $(P \ccomp Q)(g)(x \otimes q) \equiv P(g)(x) \otimes q \equiv x \otimes q$, and
    \item $(P \ccomp Q)(g)(y \otimes q) \equiv P(g)(y) \otimes q \equiv y.s.s \otimes q \equiv y \otimes Q(s.s)(q) \equiv (y \otimes q).t.t$.
\end{itemize}

In turn, $P \ccomp Q$ presents the profunctor given by
$$\semp{P \ccomp Q}(*, *) = \{[x\otimes q],\ [y\otimes q],\ [(y \otimes q).t]\}$$
$$\semp{P \ccomp Q}([f], 1) = \semp{P \ccomp Q}([g], 1): \; [x\otimes q] \mapsto [x\otimes q], \quad [y\otimes q] \mapsto [y\otimes q.t], \quad [y \otimes q.t] \mapsto [y \otimes q.t] $$
$$ \semp{P \ccomp Q}(1,[t]) : \; [x \otimes q] \mapsto [x \otimes q], \quad [y \otimes q] \mapsto [y \otimes q.t], \quad [y \otimes q.t] \mapsto [y \otimes q.t]$$
Notice that $P \ccomp Q$ is not a minimal presentation for $\semp{P\ccomp Q}$, e.g.\ the equation $(x\otimes q).t.t = (x\otimes q).t$ is redundant.
\end{Ex}

\section{Curryable Profunctor Presentations} \label{section currying}

As we have seen, finite curried profunctor presentations, unlike uncurried presentations, admit composition. However, the definition of $\ncat{Curr}(C,D)$ is more involved than that of $\ncat{UnCurr}(C,D)$. It is therefore desirable to find some subclass of uncurried profunctor presentations which behaves like the class of curried presentations.

\begin{Def} \label{def curried to uncurried}
Given a curried profunctor presentation $P:C\slashedrightarrow D$, let $\overline{P}$ denote the uncurried profunctor presentation defined as follows. Let $\Fun(\overline{P}) \coloneqq \{\overline{p} \mid c\in C, p\in\Gen(P(c)) \}$, with typing $(\overline{p} : c \to d)$ for $(p:d)\in \Gen(P(c))$. This induces a function $\overline{(-)}: \Term(P(c)) \to \CPath(\overline{P})$ for every $c \in C$, by setting $\overline{p.g} \equiv \overline{p}.g$ for every generator $p$ and $D$-path $g$. Let $\overline{P}_E \coloneqq \{ \overline{t}=\overline{t'} \mid c\in C, t=_{P(c)} t' \} + \{f.\overline{p}=\overline{P(f)(p)} \mid (f:c\to c')\in\Fun(C),\ p\in\Gen(P(c'))\}$.
Given a morphism $F:P\to Q$ of curried profunctor presentations, we define $\overline{F}:\overline{P}\to\overline{Q}$ by setting $\overline{F}(\overline{p}) :\equiv \overline{F_c(p)}$.
\end{Def}

\begin{Proprep}\label{prop curr to uncurr}
Given category presentations $C$ and $D$, the construction given in Definition \ref{def curried to uncurried} defines a functor $\overline{(-)}: \ncat{Curr}(C,D) \to \ncat{UnCurr}(C,D)$ which preserves provable equality of morphisms and makes the following diagram commute up to isomorphism
\[\begin{tikzcd}[ampersand replacement=\&]
	{\ncat{Curr}(C,D)} \&\& {\ncat{UnCurr}(C,D)} \\
	\& {\ncat{Prof}(\semp{C},\semp{D})}
	\arrow["{{\semp{-}}}"', from=1-1, to=2-2]
	\arrow["{{\semp{-}}}", from=1-3, to=2-2]
	\arrow["{{\overline{(-)}}}", from=1-1, to=1-3]
\end{tikzcd}\]
Furthermore, $\overline{(-)}$ restricts to a functor $\overline{(-)}:\ncat{FinCurr}(C,D) \to\ncat{FinUnCurr}(C,D)$.
\end{Proprep}
\begin{proof}
    We give a sketch of proof. First note that $\overline{(-)}:\Term(P(c))\to\CPath(\overline{P})$ preserves provable equality, as can be shown by an easy induction. This is the essential ingredient needed to show that whenever $F:P\to Q$ is a morphism of curried profunctor presentations, $\overline{F}$ is a morphism of uncurried profunctor presentations; and that the functor $\overline{(-)}:\ncat{Curr}(C,D)\to\ncat{UnCurr}(C,D)$ preserves provable equality of morphisms.

    The commutation up to isomorphism of the triangle also uses this same ingredient, as well as the fact that whenever $P$ is a curried profunctor presentation, $\overline{P}$ is nongenerative and conservative (see Def. \ref{def nongenerative conservative}), which will be shown in Thm. \ref{th curryable equiv to curried}.  
\end{proof}

\begin{Rem}
By Definition \ref{def curried to uncurried} and Lemma \ref{lem finite uncurried but not finite curried}, we conclude that the class of finitely curried presentable profunctors is a proper subclass of the class of finitely uncurried presentable profunctors.
\end{Rem}

\begin{toappendix}
\begin{Lemma} \label{lem provably equal to currying of cross-paths}
Let $P : C \pro D$ be a curried profunctor presentation, and let $f : c' \to c$ be a $C$-path and $t \in \Term(P(c))$. Then
\begin{equation*}
f.\overline{t} \approx_{\overline{P}} \overline{P(f)(t)}.
\end{equation*}
\end{Lemma}
\begin{proof}
We prove this by induction on $f$. If $f\equiv 1_c$, then the conclusion follows trivially.  
If the conclusion holds for a given $f:c'\to c$ and $a:c''\to c'$ is a function symbol, then $a.f.\overline{t} \approx_{\overline{P}} a.\overline{P(f)(t)}$. But $P(f)(t) \equiv q.h$ for some generator $q$ and $D$-path $h$. Thus $a.\overline{q.h} \approx_{\overline{P}} \overline{P(a)(q.h)}$ by the equations of $\overline{P}$, and since $\overline{P(a)(q.h)} \equiv \overline{P(a.f)(t)}$, the result follows.
\end{proof}
\end{toappendix}

We can begin to address our current goal by characterizing the image of the functor $\overline{(-)}:\ncat{Curr}(C,D) \to \ncat{UnCurr}(C,D)$.

\begin{Def}
Let $P : C \slashedrightarrow D$ be an uncurried profunctor presentation. By a \textbf{short left cross-path} for short, we mean a path in $P$ of the form $f.p : c \to d$, where $f: c \to c'$ is a function symbol in $C$ and $(p : c' \to d) \in \Fun(P)$. By a \textbf{right cross-path}, we mean a path in $P$ of the form $p.g : c \to d'$ where $(p : c \to d) \in \Fun(P)$ and $g: d \to d'$ is a path in $D$.
\end{Def}

\begin{Def}\label{def nongenerative conservative}
Given an uncurried profunctor presentation $P: C \xslashedrightarrow{} D$, and $c \in C$, let $P^c$ denote the $D$-instance presentation with $\Gen(P^c) = \{ (p : d) \, \mid \, p : c \to d \in \Fun(P) \}$ and $P^c_E = \{ (r=r') \in P_E \mid s(r) = c \}$. 

\end{Def}

\begin{Def}
Given an uncurried profunctor presentation $P : C \slashedrightarrow D$, we say that $P$ is \textbf{nongenerative} if for every short left cross-path $l$ in $P$, there exists a right cross-path $r$ such that $l \approx_P r$. We say that $P$ is \textbf{conservative}\footnote{The name comes from the fact that a conservative presentation $P : C \slashedrightarrow D$ is a conservative extension of each of the $P^c$, when we consider $P$ and $P^c$ as algebraic theories. See \cite[Section 4.1]{schultz2017integration}.}
if whenever $t$ and $t'$ are right cross-paths such that $t \approx_P t'$ we also have $t \approx_{P^c} t'$.
\end{Def}

\begin{Lemmarep}\label{lem nongenerative iff every cross-path equiv to right cross-path}
Let $P : C \pro D$ be a nongenerative uncurried profunctor presentation. Then for every cross-path $t$ of $P$ there exists a right cross-path $r$ such that $t \approx_P r$.
\end{Lemmarep}

\begin{proof}
Write $t \equiv f.p.g$ where $(p : c \to d)\in\Fun(P)$. Clearly the lemma holds when $f = 1_c$. Now suppose the lemma holds for a $C$-path $f:c'\to c$, and let $a:c''\to c'$ be a function symbol. By the inductive hypothesis, there exists a right cross-path $r \equiv p'.g'$, where $p'$ is a generator, such that $f.p.g \approx_P p'.g'$. Thus $a.f.p.g \approx_P a.p'.g'$. Since $P$ is nongenerative and $a.p'$ is a short left cross-path, there exists a right cross-path $r'$ such that $a.p'\approx_P r'$, so $a.f.p.g\approx_P r'.g'$.
\end{proof}

\begin{Def}
Given profunctors $\cat{P}, \cat{P}' : \cat{C} \pro \cat{D}$ and a morphism $F:\cat{P}\to \cat{P}'$, we say that $F$ is \textbf{full on cross-morphisms} if for all $c\in \cat{C}$ and $d\in \cat{D}$, the function $F_{(c,d)} : \cat{P}(c,d)\to \cat{P}'(Fc,Fd)$ is surjective.
\end{Def}

\begin{Lemmarep}
Given an uncurried profunctor presentation $P : C \slashedrightarrow D$, for every $c \in C$, there is an obvious inclusion functor\footnote{Where here we think of $\semp{P^c}$ as a profunctor $\semp{P^c} : * \pro D$.} $\iota^c : \semp{P^c} \to \semp{P}$. Thus $P$ is nongenerative if and only if for every $c \in C$, $\iota^c$ is full on cross-morphisms. Similarly $P$ is conservative if and only if for every $c \in C$, $\iota^c$ is faithful.
\end{Lemmarep}
\begin{proof}
Suppose $P$ is nongenerative. Given an element $[t]\in \semp{P}(c,d)$, by Lemma \ref{lem nongenerative iff every cross-path equiv to right cross-path} there exists a right cross-path $r$ in $P$ such that $t \approx_P r$. Since $r$ is a right cross-path starting at $c$, it is a term in $P^c$, and thus $[t] = \iota^c([r])$. Therefore $\iota^c$ is full on cross-morphisms. Conversely, suppose $\iota^c$ is full on cross-morphisms for all $c$. Then given a short left cross-path $t$ of $P$, there exists a term $r$ in $P^c$ such that $[t] = \iota^c([r])$. Since $r$ is a term in $P^c$, it is a right cross-path in $P$. This implies that $t \approx_{P} r$, and since $t$ was arbitrary, $P$ is nongenerative.

Now suppose $P$ is conservative. Given $[r], [r'] \in \semp{P^c}(d)$ with $\iota^c([r]) = \iota^c([r'])$, by conservativity we have $r \approx_{P^c} r'$, so $[r] = [r']$. Conversely, suppose $\iota^c$ is faithful. If $r$ and  $r'$ are right cross-paths in $P$ with $r\approx_P r'$, then we have $\iota^c([r]) = \iota^c([r'])$. Since $\iota^c$ is faithful, this implies that $r \approx_{P^c} r'$.

\end{proof}

\begin{Def}
We call an uncurried presentation that is both nongenerative and conservative a \textbf{curryable} profunctor presentation. We say that a morphism $F: P \to P'$ of uncurried profunctor presentations is \textbf{rightward} if it sends profunctor function symbols $p \in \Fun(P)$ to right cross-paths in $P'$. Given category presentations $C$ and $D$, let $\ncat{Crble}(C,D)$ denote the subcategory of $\ncat{UnCurr}(C,D)$ whose objects are curryable profunctor presentations and whose morphisms are rightward morphisms. Let $\ncat{FinCrble}(C,D)$ denote the full subcategory on the finite curryable profunctor presentations.
\end{Def}


Our particular definition of $\ncat{Crble}(C,D)$ is motivated by the following result.

\begin{Threp} \label{th curryable equiv to curried}
Given category presentations $C$ and $D$, the functor $\overline{(-)} : \ncat{Curr}(C,D) \to \ncat{UnCurr}(C,D)$ takes values in curryable profunctor presentations and rightward morphisms. Furthermore the corestrictions $\overline{(-)}: \ncat{Curr}(C,D) \to \ncat{Crble}(C,D)$ and $\overline{(-)}: \ncat{FinCurr}(C,D) \to \ncat{FinCrble}(C,D)$ are equivalences of categories which preserve and reflect provable equality.
\end{Threp}

\begin{proof}
Suppose that $P : C \pro D$ is a curried profunctor presentation. Then $\overline{P}$ is nongenerative by construction.  We would like to show that $\overline{P}$ is also conservative by induction using Lemma \ref{lem uncurried approx characterization}; however, to do this, we must show a statement implying conservativity which involves an arbitrary pair of provably equal $\overline{P}$-cross-paths, not just right cross-paths, as in the definition of conservativity.  Accordingly, we will show that whenever $f.\overline{t}\approx_{\overline{P}}f'.\overline{t'}$ is a provable equality of cross-paths, written such that $\overline{t}$ and $\overline{t'}$ are right cross-paths, we have $P(f)(t)\approx_{P(c)} P(f')(t')$.  Once we show this, it follows \textit{a fortiori} that whenever $\overline{t}\approx_{\overline{P}}\overline{t'}$ is a provable equality of right cross-paths, $t\approx_{P(c)} t'$, or, what is the same thing, $\overline{t}\approx_{\overline{P}^c} \overline{t'}$.

If $f.\overline{t} =_{\overline{P}} f'.\overline{t'}$, then either $f$ and $f'$ are both trivial, in which case the conclusion follows immediately, or $f$ is a function symbol and $f'$ is trivial, in which case $t'\equiv P(f)(t)$, in which case the conclusion follows by reflexivity of $\approx_{P(c)}$.  Rules 2-4 of Lemma \ref{lem uncurried approx characterization} are verified automatically since our conclusion is an equivalence relation.  The verifications for rules 5 and 7 follow immediately from the fact that $\approx_{P(c)}$ satisfies rules 5 and 6 of Lemma \ref{lem instance approx characterization}.

For rule 6 of Lemma \ref{lem uncurried approx characterization}, suppose that $f.\overline{t}\approx_{\overline{P}}f'.\overline{t'}:c\to d$ and $P(f)(t)\approx_{P(c)} P(f')(t')$, and let $(g:c'\to c)\in\Fun(C)$.  Then $P(g.f)(t)\equiv P(g)(P(f)(t))\approx_{P(c')}P(g)(P(f')(t'))\equiv P(g.f')(t')$ because $P(g)$ is a $D$-instance presentation morphism.

For rule 8 of Lemma \ref{lem uncurried approx characterization}, suppose that $f.\overline{t}\approx_{\overline{P}}f'.\overline{t'}:c\to d$ and $P(f)(t)\approx_{P(c)} P(f')(t')$, and let $g=_C g':c'\to c$.  Then $P(g.f)(t)\equiv P(g)(P(f)(t))\approx_{P(c')}P(g)(P(f')(t'))\equiv P(g.f')(t')$ because $P(g)\approx P(g')$, since $P$ is a curried presentation.

Hence $\overline{P}$ is conservative, and therefore curryable. If $F: P \to P'$ is a morphism of curried profunctor presentations, then clearly $\overline{F}$ sends profunctor function symbols to right cross-paths in $\overline{P'}$, so $\overline{F}$ is rightward, and $\overline{(-)}$ corestricts as required.

Now let us show that $\overline{(-)} : \ncat{Curr}(C,D) \to \ncat{Crble}(C,D)$ is an equivalence. Suppose that $Q:C\slashedrightarrow D$ is a curryable profunctor presentation. We must construct a curried profunctor presentation $P$ and an isomorphism $\overline{P} \cong Q$. First let us define $P:C \pro D$ sortwise by $P(c)\coloneqq Q^c$. Since $Q$ is nongenerative, for each $(f:c\to c')\in\Fun(C)$ and $p \in \Gen(P(c'))$ there exists a right cross-path $r$ such that $r \approx_Q f.p$. Using the Axiom of Choice, choose such an $r$, and set $P(f)(p) :\equiv r$. Thus we have $P(f)(p)\approx_Q f.p$, and thus holds for paths $f$ as well by induction on $f$.  Let us show that $P(f) : P(c') \to P(c)$ is a $D$-instance morphism. Given an equation $(p.g =_{P(c')} p'.g')$, we have 
\begin{equation*}
P(f)(p.g) \approx_Q f.p.g \approx_Q f.p'.g' \approx_Q P(f)(p'.g') \end{equation*}
so by conservativity, $P(f)(p.g)\approx_{Q^c} P(f)(p'.g')$ which is equivalent to $P(f)(p.g) \approx_{P(c)} P(f)(p'.g')$. Finally, given an equation $(f =_C f':c\to c')$ and a generator $p \in \Gen(P(c'))$, we have $P(f)(p) \approx_Q f.p \approx_Q f'.p \approx_Q P(f')(p)$, so by conservativity, $P(f)(p) \approx_{Q^c} P(f')(p)$, or equivalently $P(f)(p) \approx_{P(c)} P(f')(p)$, which implies that $P$ is a well defined curried profunctor presentation. Finally, we construct an isomorphism $\alpha:\overline{P} \cong Q$ by $\alpha(\overline{q}):\equiv q$, with inverse $\beta$ given by $\beta(q):\equiv \overline{q}$.  By construction $\alpha$ is a rightward morphism. Let us show that $\beta$ sends equations to provable equalities. Given an equation $f.t=_Q f'.t':c\to d$ where $t$ and $t'$ are right cross-paths, we have $P(f)(t)\approx_Q P(f')(t')$, so by conservativity $P(f)(t)\approx_{Q^c}P(f')(t')$, which is equivalent to $P(f)(t)\approx_{P(c)}P(f')(t')$, so $\beta(f.t)\equiv f.\overline{t} \approx_{\overline{P}} \overline{P(f)(t)} \approx_{\overline{P}} \overline{P(f')(t')} \approx_{\overline{P}} f'.\overline{t'}\equiv \beta(f'.t')$ by Lemma \ref{lem provably equal to currying of cross-paths}.
Thus $\overline{(-)}$ is essentially surjective.

It is easy to see that $\overline{(-)}$ is faithful; let us show it is full. Let $P, P' : C \pro D$ be curried profunctor presentations and let $G: \overline{P} \to \overline{P'}$ be a rightward morphism of curryable profunctor presentations. Thus for each $p \in \Gen(P(c))$, $G(\overline{p})$ is a right cross-path $\overline{r}$.  Thus we define $F: P \to P'$ by $F_c(p):\equiv r$.  Then for every generator we have $\overline{F}(\overline{p}) \equiv G(\overline{p})$, so $\overline{(-)}$ is full.

We have by Prop. \ref{prop curr to uncurr} that $\overline{(-)}$ preserves provable equality; to see that it also reflects provable equality, suppose that $F,G:P\to Q$ are morphisms of curried presentations with $\overline{F}\approx\overline{G}$.  Then for a generator $(p:c\to d)$ of $P$, $\overline{F_c(p)}\equiv \overline{F}(\overline{p}) \approx_{\overline{Q}} \overline{G}(\overline{p}) \equiv \overline{G_c(p)}$.  Since $\overline{Q}$ is conservative, $\overline{F_c(p)}\approx_{\overline{Q}^c} \overline{G_c(p)}$, so $F_c(p)\approx_{Q(c)}G_c(p)$ and $P\approx Q$.

It follows easily that $\overline{(-)}: \ncat{FinCurr}\to\ncat{FinCrble}$ is also an equivalence, since $\overline{(-)}$ preserves finiteness and the curried presentation $P$ constructed in our proof above of essential surjectivity is finite when $Q$ is.
\end{proof}

\begin{Cor}\label{cor finite curryable composition}
The class of profunctors which admit finite curryable profunctor presentations is closed under composition.
\end{Cor}

Conservativity is a necessary condition for Corollary \ref{cor finite curryable composition} to hold, as the following result shows.

\begin{Proprep}
There exist finite nongenerative uncurried profunctor presentations $P:C\slashedrightarrow D$ and $Q:D\slashedrightarrow E$ between finite category presentations such that $\semp{P}\odot\semp{Q}$ does not admit a finite uncurried presentation.
\end{Proprep}
\begin{proof}
Consider the monoid presentations and bimodule presentations $M\xslashedrightarrow{P}N\xslashedrightarrow{Q}O$ defined by $M\coloneqq \langle \varnothing \mid \varnothing \rangle$, $N\coloneqq \langle f \mid \varnothing \rangle$, $O\coloneqq \langle \varnothing \mid a,b \rangle$, $P(\ast)\coloneqq\langle p \mid \varnothing \rangle_N$, and $Q(\ast)\coloneqq \langle q \mid f.q=q.a,\ q=q.b\rangle_O$, as depicted in the following diagram:

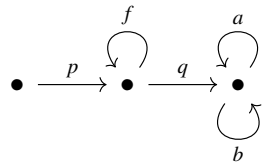
\begin{equation*}
\begin{tikzcd}
\bullet \ar[r, "p"] & \bullet \ar["f"', loop, distance=2em, in=125, out=55] \ar[r,"q"] & \bullet \ar["a"', loop, distance=2em, in=125, out=55]
\ar["b"', loop, distance=2em, in=305, out=235]
\end{tikzcd}
\end{equation*}


Both $P$ and $Q$ are nongenerative, as can be easily checked.

One can show by induction (using Lemma \ref{lem uncurried approx characterization}) that every instance of $\approx_P$ is of the form $p.f^n\approx_P p.f^n$ for some $n\in\mathbb{N}$.  Thus the equivalence classes $p_n\coloneqq [p.f^n]$ are distinct for each $n\in\mathbb{N}$ and $\semp{P}(\ast,\ast)=\{p_n\mid n\in\mathbb{N}\}$.

Similarly, induction shows that every instance of $\approx_Q$ is of the form $f^n.q.w\approx_Q f^{n'}.q.w'$ satisfying $n+\#_a(w)=n'+\#_a(w')$, where $\#_a(w)$ is the number of times that $a$ appears in the path $w$.  Thus the equivalence classes $q_n\coloneqq [q.a^n]$ are distinct for $n\in\mathbb{N}$ and $\semp{Q}(\ast,\ast)=\{q_n\mid n\in\mathbb{N}\}$, with $\semp{Q}(f,1)(q_n)=q_{n+1}$, $\semp{Q}(1,a)(q_n)=q_{n+1}$, and $\semp{Q}(1,b)(q_n)=q_n$. 

Now consider the composite profunctor $\semp{P} \odot \semp{Q}: \semp{C} \slashedrightarrow \semp{E}$.  The relation $\simeq_0$ as in \eqref{eq description of profunctor composition hom} is then defined by $(p_m,q_{n+o})\simeq_0 (p_{m+o},q_n)$ for $m,n,o\geq 0$.
Thus $(p_m,q_n)\simeq (p_{m'}, q_{n'})$ whenever $m+n=m'+n'$, and moreover, $\{((p_m,q_n),\ (p_{m'}, q_{n'})) \mid m+n=m'+n' \}$ is an equivalence relation so it is exactly $\simeq$.  Hence the equivalence classes $pq_n\coloneqq \langle p_0,q_n\rangle$ are distinct for $n\in\mathbb{N}$ and $(\semp{P} \odot \semp{Q})(\ast,\ast)=\{pq_n \mid n\in\mathbb{N} \}$, with $(\semp{P} \odot \semp{Q})(1,a)(pq_n)=pq_{n+1}$ and $(\semp{P} \odot \semp{Q})(1,b)(pq_n)=pq_n$.

Now suppose that $R$ is a finite uncurried presentation and $K:\semp{R}\cong \semp{P} \odot \semp{Q}$ is an isomorphism. Define the \textbf{order} $o(t)$ of a cross-path $t$ of $R$ as the unique $n$ such that $K([t])=pq_n$.  Note that provably equal cross-paths of $R$ have identical order.  Moreover, since $K$ is an isomorphism, the function $o:\CPath(R)\to\mathbb{N}$ is surjective.

Since $R_E$ is finite, let $N\coloneqq\max\{o(t) \mid t=_R t'\}$ (we let $N\coloneqq 0$ if $R_E$ is empty).  It is not hard to show by induction that whenever $t\approx_R t'$, either $o(t)\leq N$ or there exist $s,s'\in\CPath(R)$ and $g\in\Fun(O)$ such that $t\equiv s.g$ and $t'\equiv s'.g$.  Let $t\in\CPath(R)$ be the shortest cross-path such that $o(t)=N+1$.  Then we must have $[t.b]=\semp{R}(1,b)([t])=[t]$, so $t.b\approx_R t$.  Since $o(t)>N$, there must be $s\in\CPath(R)$ such that $t=_R s.b$.  But then $s$ is a shorter cross-path than $t$ with $o(s)=N+1$, contradicting the definition of $t$.
\end{proof}

\begin{Rem}
We can revisit Example \ref{ex curried example} replacing curried profunctor presentations $P$ and $Q$ for their uncurried counterparts $\overline{P}$ and $\overline{Q}$.
We can then visualize the overall situation as follows:

\vspace*{-0.8cm}
\begin{multicols}{3}

\vspace*{\stretch{0}}
\begin{equation*}
{\footnotesize
\begin{aligned}
    f.g&=g.f\phantom{\overline{\otimes}} \\
    f.\overline{y}&=\overline{y}.s\phantom{\overline{\otimes}} \\
    f.\overline{x}&=g.\overline{x}=\overline{x}\phantom{\overline{\otimes}} \\
    \overline{x}.s&=\overline{x}\phantom{\overline{\otimes}} \\
    g.\overline{y}&=\overline{y}.s.s\phantom{\overline{\otimes}} \\
    \overline{q}.t.t &= \overline{q}.t\phantom{\overline{\otimes}} \\
    s.\overline{q} &= \overline{q}.t\phantom{\overline{\otimes}}
\end{aligned}
}
\end{equation*}
\vspace*{\stretch{1}}

\columnbreak

\begin{equation*}
\begin{tikzcd}[row sep = 30pt]
\bullet \arrow["g"', loop, distance=2em, in=215, out=145] \arrow["f"', loop, distance=2em, in=125, out=55] \arrow[r, "{\overline{x}}", bend left] \arrow[r, "{\overline{y}}"', bend right] & \bullet \arrow[r, "{\overline{q}}"] \arrow["s"', loop, distance=2em, in=125, out=55] & \bullet \arrow["t"', loop, distance=2em, in=125, out=55] \\
\bullet \arrow["g"', loop, distance=2em, in=215, out=145] \arrow["f"', loop, distance=2em, in=125, out=55] \arrow[rr, "{\overline{x \otimes q}}", bend left=20] \arrow[rr, "{\overline{y \otimes q}}"', bend right=20] && \bullet \arrow["t"', loop, distance=2em, in=125, out=55]
\end{tikzcd}
\end{equation*}

\columnbreak

\vspace*{\stretch{2}}
\begin{equation*}
{\footnotesize
\begin{aligned}
    (\overline{x\otimes q}).t.t &= (\overline{x\otimes q}).t \\
    (\overline{y\otimes q}).t.t &= (\overline{y\otimes q}).t \\
    (\overline{x\otimes q}).t &= \overline{x\otimes q} \\
    f.(\overline{x\otimes q}) &= \overline{x\otimes q} \\
    f.(\overline{y\otimes q}) &= (\overline{y\otimes q}).t \\
    g.(\overline{x\otimes q}) &= \overline{x\otimes q} \\
    g.(\overline{y\otimes q}) &= (\overline{y\otimes q}).t
\end{aligned}    
}
\end{equation*}
\vspace*{\stretch{1}}

\end{multicols}

\vspace*{0.1cm}

The diagram on top, together with the equations on the left, represents the situation before computing the composite, i.e. $\overline{P}$ and $\overline{Q}$ drawn together. Meanwhile, the bottom diagram together with the equations on the right represents $\overline{P \ccomp Q}$. This perspective on composition of profunctor presentations, using curryable presentations, has the benefit of being more easily visualized.

\end{Rem}

\begin{Rem} \label{rem other interpretations of main theorem 2}
We have chosen to present a result on a category of curryable profunctor presentations from $C$ to $D$ which is \emph{equivalent} to $\ncat{Curr}(C,D)$. However, other options are possible, depending on the motivation.
    
From a theoretical perspective, we may find the restriction to rightward morphisms undesirable, since it makes $\ncat{Crble}(C,D)$ a non-full subcategory of $\ncat{Curr}(C,D)$. We can solve this by relaxing our notion of equivalence. If we now consider $\ncat{Crble}(C,D)$ to be the full subcategory of $\ncat{UnCurr}(C,D)$ spanned by curryable presentations, then $\overline{(-)}:\ncat{Curr}(C,D)\to\ncat{Crble}(C,D)$ is no longer an equivalence, but is a \emph{biequivalence} between bicategories (considering provable equalities as 2-morphisms as in Remark \ref{rem bicategories of uncurried profunctor presentations}) which restricts to a biequivalence between the corresponding full sub-bicategories of finite presentations.\footnote{Our original functor $\overline{(-)}: \ncat{Curr}(C,D) \to \ncat{Uncurr}(C,D)$ is already a biequivalence in this way, but it does not restrict to a biequivalence between the full sub-bicategories of finite presentations.}

Alternatively, from the point of view of implementation we may be interested in being able to replace curried presentations by uncurried presentations through a bijective correspondence. In this case we could then strengthen nongenerativity to require that all equations are either between right cross-paths or between a short left cross-path and a right cross-path, and every short left cross-path appears in some equation exactly once. Using this notion of nongenerativity, the functor $\overline{(-)}:\ncat{Curr}(C,D)\to\ncat{Crble}(C,D)$ actually becomes an isomorphism of categories, which again restricts to an isomorphism between the corresponding full subcategories of finite presentations.
\end{Rem}

\section{The Double Category of Curried Profunctor Presentations}\label{section double categories}

\begin{toappendix}
We include definitions of double categories and double functors for ease of reference and to fix notation for the proofs that follow.

\begin{Not}
In what follows, given morphisms $f: A \to B$ and $g: C \to B$, let $A \times_{f, \, B, \, g} C$ denote the pullback of the cospan $A\rightarrow B \leftarrow C$, with the canonical projections $\pi:A \times_{f, \, B, \, g} C \to A$ and $\pi': A \times_{f, \, B, \, g} C \to C$.
\end{Not}

\begin{Def} \label{def double category}
A \textbf{double category}\footnote{By a double category, we mean what is sometimes called a weak or pseudo double category in the literature \cite{grandis2019higher}. Namely we consider horizontal arrow composition to be weakly unital and associative, in contrast to strict double categories} $\mathbb{D}$ consists of
\begin{itemize}
    \item A category $\mathbb{D}_0$, which we refer to as the \textbf{vertical category} of $\mathbb{D}$.
    \item A category $\mathbb{D}_1$, equipped with two functors $L, R: \mathbb{D}_1 \to \mathbb{D}_0$, called the \textbf{left
and right frame functors}. Given an object $\cat{P} \in \mathbb{D}_1$ with $\cat{C} = L(\cat{P})$ and $\cat{D} = R(\cat{P})$, we say that $\cat{P}$ is a \textbf{proarrow} (or \textbf{horizontal arrow}) from $\cat{C}$ to $\cat{D}$ and write $\cat{P} : \cat{C} \slashedrightarrow \cat{D}$. We call morphisms $\phi : \cat{P} \to \cat{P}'$ in $\mathbb{D}_1$ \textbf{2-cells}. We say a $2$-cell $\phi$ is \textbf{globular} if $L(\phi)$ and $R(\phi)$ are identity morphisms in $\mathbb{D}_0$.
\item A functor $U : \mathbb{D}_0 \to \mathbb{D}_1$, which we call the \textbf{unit}, such that $L \circ U = R \circ U = 1_{\mathbb{D}_0}$.
\item A functor $\odot : \mathbb{D}_1 \times_{R, \, \mathbb{D}_0, \, L} \mathbb{D}_1 \to \mathbb{D}_1$ called \textbf{horizontal composition}, such that $L\circ\odot=L\circ\pi$ and $R\circ\odot=R\circ \pi'$, which is weakly
associative and weakly unital in the sense that there are coherent, globular left and right unitor and associator isomorphisms \cite[Definition 12.3.7]{johnson2020} denoted $\lambda$, $\rho$, and $\alpha$, respectively.

\end{itemize}
Given a double category $\mathbb{D}$ and objects $c, d \in \mathbb{D}_0$, let $\mathbb{D}(c,d)$ denote the subcategory of $\mathbb{D}_1$ whose objects are proarrows $p \in \mathbb{D}_1$ such that $L(p) = c$ and $R(p) = d$, and whose morphisms are globular $2$-cells $\phi : p \to p'$.
\end{Def}

\begin{Def}
A strong double functor $\mathbb{F}:\mathbb{D}\to\mathbb{E}$ of double categories consists of
\begin{itemize}
    \item A functor $\mathbb{F}_0:\mathbb{D}_0\to\mathbb{E}_0$.
    \item A functor $\mathbb{F}_1:\mathbb{D}_1\to\mathbb{E}_1$ such that $L\circ \mathbb{F}_1 = \mathbb{F}_0\circ L$ and $R\circ \mathbb{F}_1 = \mathbb{F}_0\circ R$.
    \item Globular isomorphisms $\mu:\mathbb{F}_1\circ \odot \cong \odot \circ \mathbb{F}_1\times_{\mathbb{F}_0}\mathbb{F}_1$ and $\eta:\mathbb{F}_1\circ U\cong U\circ\mathbb{F}_0$ which are coherent with respect to the left and right unitors and associators of $\mathbb{D}$ and $\mathbb{E}$ (see \cite{ferreira2006pseudocategories} for the full definition).
    \end{itemize}
\end{Def}
In this paper we will only consider double functors that are strong, and thus omit the word ``strong.''

\begin{Def} \label{double congruence}
Let $\mathbb{D}$ be a double category, and let $\sim_i$ be a congruence on $\mathbb{D}_i$ for $i=0,1$.  If $L$, $R$, $U$, and $\odot$ preserve these congruences, we say that $\sim=(\sim_0,\sim_1)$ is a double congruence on $\mathbb{D}$.  Then we can define a double category $\mathbb{D}/\sim$, by $\mathbb{D}_i/\sim_i$ for $i=0,1$ and the functors induced by $L$, $R$, $U$ and $\odot$.
\end{Def}

\begin{Lemma}\label{double congruence universal}
Given a double category $\mathbb{D}$ with a double congruence $\sim$, every double functor $F:\mathbb{D}\to\mathbb{D'}$ factors uniquely through $\mathbb{D}/{\sim}$ to induce a functor $(\mathbb{D}/{\sim})\to\mathbb{D'}$.  
\end{Lemma}

\begin{Def}[Non-globular morphisms]
Let $\ncat{Prof} := \ncat{Cat}/\ncat{2}$ and let $\ncat{Inst}$ be the full subcategory of $\ncat{Prof}$ on instances. Now define the category $\ncat{InstPr}$ as follows: its objects are $C$-instance presentations for some $C \in \ncat{CatPr}$, and a morphism $F$ from $I \in C\hyp\ncat{InstPr}$ to $J \in D\hyp\ncat{InstPr}$ is a morphism of category presentations $|I|\to |J|$ which sends $\ast\mapsto\ast$ and sorts of $C$ to sorts of $D$. In other words, it is given by a morphism of category presentations $F_1: C \to D$ and a function $F_{01}:\Gen(I)\to\Term(J)$ such that $p : c$ implies $F_{01}(p): F_1(c)$ and for all equations $p.f =_I p'.f'$, $F_{01}(p).F_1(f) \approx_J F_{01}(p').F_1(f')$. Let $\ncat{InstPr}_\approx$ be defined by quotienting by $\approx$, which holds between $F$ and $G$ if and only if $F_1 \approx G_1$ and $F_{01}(p) \approx_J G_{01}(p)$ for all $p \in \Gen(I)$.

Finally, define the category $\ncat{Curr}$ as follows: its objects are curried profunctor presentations, and a morphism $F$ from $P : C \pro D$ to $P' : C' \pro D'$ consists of category presentation morphisms $F_0:C\to C'$ and $F_1:D\to D'$, along with a collection of morphisms $F_{c}:P(c)\to P'(F_0(c))$ in $\ncat{InstPr}$ with $(F_{c})_1=F_1$, indexed by the sorts $c$ of $C$, such that for each function symbol $f: c \to c'$ in $C$, the following diagram commutes up to provable equality:
    \begin{equation} \label{eq diagram of nonglobular morphism of curried presentations}
    \begin{tikzcd}[ampersand replacement=\&]
        {P(c')} \& {P'(F_0(c'))} \\
        {P(c)} \& {P'(F_0(c))}
        \arrow["{F_{c'}}", from=1-1, to=1-2]
        \arrow[""{name=0, anchor=center, inner sep=0}, "{P(f)}"', from=1-1, to=2-1]
        \arrow[""{name=1, anchor=center, inner sep=0}, "{P'(F_0(f))}", from=1-2, to=2-2]
        \arrow["{F_{c}}"', from=2-1, to=2-2]
        \arrow["\approx"{description, pos=0.4}, draw=none, from=0, to=1]
    \end{tikzcd}
    \end{equation}
    Notice that \eqref{eq diagram of nonglobular morphism of curried presentations} still commutes when $f$ is a path, by induction on $f$; thus these morphisms can be composed in the natural way.
    Let $\ncat{FinCurr}$ denote the full subcategory of finite curried profunctor presentations, i.e. $P$ such that $C$, $D$, and $P(c)$ for each $c$ are finite.
    We say that two curried profunctor presentation morphisms $F,G:P\to P'$ are \textbf{provably equal}, and write $F\approx G$, if $F_0 \approx G_0$, $F_1 \approx G_1$, and $F_c\approx G_c:P(c)\to P'(F_0(c))=P'(G_0(c))$ for each $c\in C$.
\end{Def}

\begin{Def}
Let $L,R:\ncat{Prof}\to\ncat{Cat}$ be the functors defined objectwise by sending a profunctor $\cat{P} : \cat{C}\slashedrightarrow \cat{D}$ to the categories $\cat{C}$ and $\cat{D}$ respectively. Then the operation $\odot$ of profunctor composition in Definition \ref{def profunctor composition} extends to a functor 
${\odot : \ncat{Prof}\times_{R, \, \ncat{Cat}, \, L} \ncat{Prof} \to \ncat{Prof}.}$
Let $U: \ncat{Cat} \to \ncat{Prof}$ be defined on objects by sending $\cat{C}$ to the profunctor $\text{pr}_2: \cat{C} \times \ncat{2} \to \ncat{2}$ and on morphisms by sending a functor $F: \cat{C} \to \cat{D}$ to the morphism of profunctors $F \times \ncat{2}: \cat{C} \times \ncat{2} \to \cat{D} \times \ncat{2}$.
Together, the functors $L, R, U$ and $\odot$ form a double category $\bbprof$ such that $\bbprof_0 \coloneqq \ncat{Cat}$ and $\bbprof_1 \coloneqq \ncat{Prof}$. In other words, the objects are (small) categories, vertical morphisms are functors, horizontal morphisms are profunctors, with horizontal composition given by profunctor composition and the $2$-cells are morphisms of profunctors. Note that $\bbprof(\cat{C}, \cat{D})$ is the subcategory of $\ncat{Prof}$ whose objects are profunctors $\cat{P}: \cat{C} \slashedrightarrow \cat{D}$, and whose morphisms $F : \cat{P} \to \cat{P}'$ restrict to the identity on $\cat{C}$ and $\cat{D}$. See \cite[Example 2.12]{schultz2017algebraic} and \cite[Remark 5.1.8]{loregian2021co} for more details.
\end{Def}

\end{toappendix}

A natural categorical structure to think about profunctors is that of double categories. By a double category, we mean what is sometimes called a weak or pseudo double category in the literature \cite{ferreira2006pseudocategories,fiore2007pseudo}. Namely we consider horizontal arrow composition to be weakly unital and associative, in contrast to strict double categories \cite[Explanation 12.3.1]{johnson2020}. The double category $\bbprof$ packages together the morphisms between profunctors and the composition operation of profunctors. Analogously, curried profunctor presentations can also be arranged in this way together with $\circledast$ into a double category $\bbcurr$.



\begin{toappendix}
\begin{Def}
Let $L_{\ncat{Curr}},R_{\ncat{Curr}}:\ncat{Curr}\to\ncat{CatPr}$ denote the left and right frame functors, sending a curried profunctor presentation $P : C \slashedrightarrow D$ to the category presentations $C$ and $D$ respectively.
\end{Def}

\begin{Lemma}
    The construction $P \mapsto \semp{P}$ of Definition \ref{def semantics of curried profunctor presentations, globular case} extends to a functor $\semp{-}: \ncat{Curr} \to \ncat{Prof}$.
\end{Lemma}
\begin{proof}
    Given curried profunctor presentations $P: C \pro D$ and $P' : C' \pro D'$ and a morphism $\phi: P \to P'$ between them, let $\semp{\phi}: \semp{P} \to \semp{P'}$ be defined as the following functor of categories over $\ncat{2}$. On objects and morphisms of $\semp{C}$ it acts as $\semp{\phi_0}: \semp{C} \to \semp{C'}$, while on objects and morphisms of $\semp{D}$ it acts as $\semp{\phi_1}: \semp{D} \to \semp{D'}$. Finally, given any cross-morphism
    $[t] \in \semp{P(c)}(d) = \semp{P}(c,d)$,
    let $\semp{\phi}([t]) \coloneq [\phi_{c}(t)]$. It is straightforward to see that this gives a well defined functor which commutes with the projection functors into $\ncat{2}$. Functoriality is also immediate from the definition.
\end{proof}


\begin{Lemma}\label{lem composition of curried presentations is functorial non globular case}
The assignment $(P, Q) \mapsto P \ccomp Q$ defines a functor $\ccomp : \ncat{Curr}\times_{R_{\ncat{Curr}}, \, \ncat{CatPr}, \, L_{\ncat{Curr}}}\ncat{Curr} \to \ncat{Curr}.$
\end{Lemma}
\begin{proof}
This is just a generalization of Lemma \ref{lem composition of curried presentations is functorial} to non-globular morphisms of curried profunctor presentations. Accordingly, the proof is completely analogous. For concreteness, we give the action of $\ccomp$ on non-globular morphisms. 

Given curried profunctor presentations $P : C \slashedrightarrow D$, $P' : C' \slashedrightarrow D'$, $Q: D \slashedrightarrow E$ and $Q' : D' \slashedrightarrow E'$, and given non-globular morphisms $\phi: P \to P'$ and $\psi: Q \to Q'$ such that $R_{\ncat{Curr}}(\phi)=L_{\ncat{Curr}}(\psi)$, let $\phi\ccomp \psi: P\ccomp Q\to P'\ccomp Q'$ be defined by $(\phi\ccomp \psi)_0\coloneqq\phi_0:C\to C'$, and $(\phi\ccomp \psi)_1\coloneqq\psi_1:E\to E'$, and for each $c\in \Sort(C)$, let $(\phi \ccomp \psi)_c:(P\ccomp Q)(c)\to (P'\ccomp Q')(\phi_0(c))$ be defined as in Lemma \ref{lem composition of curried presentations is functorial} by sending generators $(p \otimes q)$ to $(\phi_c(p) \otimes \psi_d(q))$, where $d=t(p)$. Using this definition, the proof follows the same steps as those of Lemma \ref{lem composition of curried presentations is functorial}.

\end{proof}
\end{toappendix}

\begin{Threp}
There exists a double category $\bbcurr$ such that the objects are category presentations and the horizontal morphisms are curried profunctor presentations.
\end{Threp}
\begin{proof}
Let us show that the functors $L_{\ncat{Curr}}$, $R_{\ncat{Curr}}$, and $\circledast$ form a double category $\bbcurr$ whose objects are category presentations, with vertical morphisms given by morphisms of category presentations, horizontal morphisms given by curried profunctor presentations and $2$-cells given by morphisms of curried profunctor presentations.

We first define a unit functor $U:\ncat{CatPr}\to\ncat{Curr}$. For a category presentation $C$ and a sort $c\in C$, we define $U(C)(c)$ to be the $C$-instance presentation with one generator $(x_c:c)$ and no equations. Given a function symbol $(f:c\to c')$ in $C$, we define $U(C)(f): U(C)(c')\to U(C)(c)$ by $(x_{c'}\mapsto x_c.f)$. Given a category presentation morphism $F:C\to C'$, we define $U(F):U(C)\to U(C')$ by $U(F)_0=U(F)_1=F$ and, for each $c\in C$, we define $U(F)_c:U(C)(c)\to U(C')(Fc)$ by $(x_c\mapsto x_{Fc})$. Clearly, we have $L\circ U = R\circ U = 1_{\ncat{CatPr}}$.

Next, we define the globular associator and unitor isomorphisms $\alpha^{P,Q,R}:(P\circledast Q)\circledast R \cong P\circledast (Q\circledast R)$, $\lambda^P: I(L(P))\circledast P \cong P$, and $\rho^P: P\circledast I(R(P))\cong P$, natural in all variables.

For a generator $((p\otimes q)\otimes r)\in\Gen(((P\circledast Q)\circledast R)(c))$, we define $\alpha^{P,Q,R}_c((p\otimes q)\otimes r) :\equiv p\otimes (q\otimes r)$. Then for terms $(s:d)\in \Term(P(c))$, $(t:e)\in\Term(Q(d))$, and $(u:f)\in\Term(R(e))$, let $s\equiv p.g$, $Q(g)(t)\equiv q.h$, and $R(h)(u)\equiv r.k$. Then $\alpha^{P,Q,R}_c$ maps $(s\otimes t)\otimes u\equiv (p\otimes Q(g)(t))\otimes u \equiv (p\otimes q).h\otimes u \equiv ((p\otimes q)\otimes r).k$ to $(p\otimes (q\otimes r)).k$, which is easily shown to be syntactically equal to $s\otimes (t\otimes u)$. Using this fact we can easily verify that $\alpha^{P,Q,R}$ is a curried profunctor presentation morphism, and it is also clearly invertible and globular. This fact also makes naturality in $P$, $Q$, and $R$ straightforward to show.

Next, for a generator $p\in\Gen(P(c))$, we define $\lambda^P_c(x_c\otimes p):\equiv p$. Clearly $\lambda^P_c$ is an invertible $C$-instance morphism. Given a $C$-function symbol $f:c\to c'$, we have $\lambda^P_{c}(U(P)(f)(x_{c'})\otimes p) \equiv \lambda^P_{c}(x_c\otimes P(f)(p)) \equiv P(f)(p) \equiv P(f)(\lambda^P_{c'}(x_c\otimes p))$, so $\lambda^P$ is an isomorphism, and it is easy to see that it is natural in $P$ and globular.

For a generator $(p:d)\in\Gen(P(c))$, we define $\rho^P_c(p\otimes x_d):\equiv p$. For a $D$-path $g:d\to d'$, we have $\rho^P_c(p.g\otimes x_{d'})\equiv \rho^P_c(p\otimes x_{d}.g)\equiv p.g$. Just as above, $\rho^P$ is an isomorphism, natural in $P$, and globular.

The pentagon and triangle identities are trivially satisfied.
\end{proof}

The double category $\bbcurr$ organizes the composition of curried profunctor presentations and their transformations in a way that mimics the structure of $\bbprof$. The isomorphism of Theorem \ref{th curried composition semantics} is now packaged into a double functor between these double categories.

\begin{Threp}\label{thm double functor from curried to profunctors}
The construction $\semp{-}$ on curried profunctor presentations defines a double functor $\semp{-}: \bbcurr \to \bbprof$.
\end{Threp}
\begin{proof}
Let $\semp{-}_0: \bbcurr_0 \to \bbprof_0$ and $\semp{-}_1: \bbcurr_1 \to \bbprof_1$ be given by the semantic functors $\semp{-}:\ncat{CatPr}\to\ncat{Cat}$ and $\semp{-}:\ncat{Curr}\to\ncat{Prof}$. Notice that indeed $L \circ \semp{-}_1 = \semp{-}_0 \circ L_{\ncat{Curr}}$ and $R \circ \semp{-}_1 = \semp{-}_0 \circ R_{\ncat{Curr}}$ as required. The desired globular isomorphism $\mu: \odot \circ \semp{-}_1 \times_{\semp{-}_0} \semp{-}_1 \cong \semp{-}_1 \circ \ccomp$ is essentially given by the isomorphism of Theorem \ref{th curried composition semantics}, generalized to non-globular morphisms of curried profunctor presentations. We now describe this more precisely.

We are looking for a natural isomorphism $\mu: \semp{-} \odot \semp{=} \Rightarrow \semp{- \ccomp =}$ between functors of the type $\ncat{Curr} \times_{R_{\ncat{Curr}}, \ncat{CatPr}, L_{\ncat{Curr}}} {\ncat{Curr}} \to \ncat{Prof}$. This means that we must give a component isomorphism $\mu^{P,Q}: \semp{P} \odot \semp{Q} \Rightarrow \semp{P\ccomp Q}$ for each pair of composable curried profunctor presentations $P$ and $Q$. This is already provided by Theorem \ref{th curried composition semantics}. The only difference, then, is that in order to establish naturality we must now consider curried profunctor presentations $P: C \pro D$, $P': C' \pro D'$, $Q: D \pro E$ and $Q': D'\pro E'$, together with arbitrary morphisms $\phi: P \to P'$ and $\psi: Q \to Q'$, and show that $\semp{\phi \ccomp \psi} \circ \mu^{P,Q} = \mu^{P',Q'} \circ (\semp{\phi}\odot\semp{\psi})$ using the functorial action of $\ccomp$ on non-globular morphisms given in Lemma \ref{lem composition of curried presentations is functorial non globular case}. This is equally immediate as the verification in the globular case. To clarify the role of the non-globular morphisms we repeat the verification expanding, without loss of generality, $s \equiv p$ and $t \equiv q.h$ for some $(p:d) \in \Gen(P(c))$, $(q:e') \in \Gen(Q(d))$ and $E$-path $h:e'\to e$. We then have that
\begin{align*}
(\semp{\phi \ccomp \psi} \circ \mu^{P,Q}) (\langle [p], [q.h]\rangle) &= \semp{\phi \ccomp \psi} ([p \otimes q.h]) \\
&= [\phi_c(p) \otimes \psi_d(q).\psi_1(h)] \\
&= \mu^{P',Q'}(\langle [\phi_c(p)], [\psi_d(q).\psi_1(h)]\rangle) \\
&= (\mu^{P',Q'} \circ (\semp{\phi}\odot\semp{\psi})) (\langle[p],[q.h]\rangle).
\end{align*}


Now we define a globular isomorphism $\varepsilon^C:\semp{U(C)}\to \Hom_{\semp{C}}$ as follows. Given a $C$-path $f:c\to c'$, let $\varepsilon^C_c([x_c.f])=[f]$. To see that this is well-defined, we must show that for $f,f':c\to c'$, 
\begin{equation}
        \left(x_c.f\approx_{U(C)(c)}x_c.f'\right) \Longleftrightarrow \left(f\approx_C f'\right)
\end{equation}

    For the $(\Rightarrow)$ direction, we observe that $\approx_C$ is an equivalence relation, includes $(f,f')$ for every $(x_c.f=x_c.f')\in U(C)(c)_E=\varnothing$, includes $C$-equations, and is closed under right-composition with a function symbol. For the $(\Leftarrow)$ direction, we simply observe that $\approx_{|U(C)(c)|}$ specializes to both $\approx_C$ and $\approx_{U(C)(c)}$ and respects left-composition with $x_c$.

    It is easily observed that $\varepsilon_c^C$ is a globular morphism of profunctors, and that it is natural in $c$ and $C$.

    Finally, we need to check that $\mu$ and $\epsilon$ respect $\alpha$, $\lambda$, and $\rho$. However, this is straightforward.
\end{proof}

Analogously to the case of $\ncat{Curr}(C,D)$, we can quotient $\bbcurr$ by an appropriate congruence of double categories induced by the provable equality relation.
\begin{Proprep}
The congruences $\approx$ on $\ncat{Curr}$ and $\ncat{CatPres}$ comprise a double congruence on $\bbcurr$.  Define $\bbcurr_\approx$ by quotienting $\bbcurr$ by $\approx$, as in Def. \ref{double congruence}. The semantics double functor $\semp{-}$ induces a double equivalence $\semp{-}:\bbcurr_{\approx} \to \bbprof$.
\end{Proprep}
\begin{proof}
The semantics double functor $\semp{-}$ is invariant with respect to provable equality, so by Lemma \ref{double congruence universal} there is an induced double functor $\semp{-}:\bbcurr/{\approx} \to \bbprof$. To show that this strong double functor is an equivalence, \cite[Theorem 7.8]{shulman2009framed} tells us that we can show that it is full, faithful, and essentially surjective in the sense of \cite[Def. 3.3, 3.6]{shulman2009framed}. This is straightforward.
\end{proof}

We can also express the fact that composition of curried profunctor presentations preserves finiteness as the existence of a double subcategory $\mathbb{F}\ncat{inCurr}\subseteq \bbcurr$.

In future work we plan to study syntactic objects using double categories equipped with double congruences as extra structure in a way analogous to Remark \ref{rem bicategories of uncurried profunctor presentations}, but in this case leading to a three-dimensional structure. We believe that this could constitute a general framework for reasoning about presentations of proarrows in double categories. 


\printbibliography

\end{document}